\documentclass{amsart}
\usepackage{amssymb}
\usepackage{color}
\usepackage{verbatim} 
\usepackage{amsmath,amssymb}

\newtheorem{The}{Theorem}[section]
\newtheorem{Pro}[The]{Proposition}
\newtheorem{lem}[The]{Lemma}

\newtheorem{Def}[The]{Definition}
\newtheorem{Rem}[The]{Remark}
\newtheorem{Examp}[The]{Example}

\usepackage{verbatim} 
\def\a1s{a_1,\cdots, a_s}
\def\a{\alpha}

\def\andd{\quad\hbox{and}\quad}

\def\b{\beta}

\def\bl4{B_{\ell\geq4}}

\def\d{\delta}

\def\bbbf{\mathbb{F}}

\def\bbbk{\mathbb{K}}
\def\lam{\lambda}

\def\LL{\mathcal{L}}

\def\ep{\epsilon}
\def\fm{(\cdot,\cdot)}

\def\bbbq{\mathbb{Q}}

\def\bbbn{\mathbb{N}}

\def\rd{\dot{R}}

\def\rds{\dot{R}_{sh}}
\def\rdl{\dot{R}_{lg}}
\def\rde{\dot{R}_{ex}}

\def\1k{\frac{1}{k}}
\def\op{\oplus}
\def\ot{\otimes}

\def\la{\langle}
\def\ra{\rangle}

\def\sub{\subseteq}
\def\sg{\sigma}

\def\rcross{R^{\times}}

\def\pf{\noindent{\bf Proof. }}

\def\T{{\mathcal T}}

\def\u{{\mathcal U}}

\def\v{{\mathcal V}}

\def\w{{\mathcal W}}

\def\bbbz{{\mathbb Z}}

\def\1il{1\leq i\leq\ell}

\def\rre{R_{re}}
\def\rim{R_{ns}}

\begin{document}

\title{Extended Affine Root Supersystems}
\thanks{Address: Department of Mathematics, University of Isfahan, Isfahan, Iran,
P.O.Box 81745-163, and School of Mathematics, Institute for Research in
Fundamental
Sciences (IPM), P.O. Box: 19395-5746, Tehran, Iran.\\
Email addresse: ma.yousofzadeh@sci.ui.ac.ir \& ma.yousofzadeh@ipm.ir.}
\maketitle
\centerline{Malihe Yousofzadeh}
\vspace{2cm}

\begin{abstract}
 The interaction of a Lie algebra $\LL,$ having a weight space decomposition with respect to a nonzero toral subalgebra, with its corresponding root system forms  a powerful tool in the study  of the  structure of $\LL.$ This, in particular, suggests a  systematic study of the root system apart from its connection with the Lie algebra. Although there have been a lot of  researches  in this regard  on Lie algebra level, such an approach has not been considered on Lie superalgebra level.   In this work, we  introduce  and study extended affine  root supersystems which are  a generalization of both  affine  reflection  systems and  locally finite root supersystems. Extended affine root supersystems appear as the root systems of the  super version of extended affine Lie algebras and invariant affine reflection algebras including affine Lie superalgebras. This work  provides a framework to  study the structure of  this kind of Lie superalgebras  refereed to as extended affine Lie superalgebras.
\end{abstract}

\setcounter{section}{-1}

\section{Introduction}
 Lie algebras having a weight space decomposition with respect to a nonzero abelian subalgebra, called  a toral subalgebra, form a vast class of Lie algebras. Locally finite split simple Lie algebras \cite{NS}, extended affine Lie algebras \cite{AABGP}, toral type extended affine Lie algebras \cite{AKY}, locally extended affine Lie algebras \cite{MY} and  invariant affine reflection algebras \cite{N1} are examples of such Lie algebras.  We can attach to such a Lie algebra, a subset of the dual space of its toral subalgebra called the  root system. The interaction of such a Lie algebra with its root system offers an approach to study the structure of the Lie algebra via its root system. This in turn provokes  a  systematic study of the  root system apart form its connection with the Lie algebra; see \cite{AABGP}, \cite{LN}, \cite{yos2} and \cite{N1}. Although since 1977, when the concept of Lie superalgebras was introduced  \cite{K1}, there has been  a significant number of researches on Lie superalgebras, the mentioned approach on   Lie superalgebra level has not been considered in general. The first step towards such an  approach is offering  an abstract definition of the root system of a Lie superalgebra. In  1996, V. Serganova \cite{serg} introduced the notion of  generalized root systems as a  generalization of finite root systems; see also \cite{HY}.
The main difference between generalized root systems and  finite root
systems is the existence of nonzero self-orthogonal roots.  {Serganova classified irreducible generalized root systems and  showed  that  such  root systems are root systems
of finite dimensional basic classical simple Lie  superalgebras  \cite{K1} except for type $A(1,1).$ She also  gave  two alternative definitions for generalized root systems. In a generalized root system  for two  self-orthogonal roots  which are not orthogonal, either their summation or their subtraction (and not both) is again a root while according to  the first alternative definition both summation and subtraction of two  self-orthogonal roots  which are not orthogonal, can be roots; this in particular allows  to obtain type $A(1,1)$ as well.
}
In this work, we introduce extended affine root supersystems and systematically study them. Roughly speaking, a spanning set    $R$ of a  nontrivial vector space over a field $\bbbf$ of characteristic zero, equipped with a  symmetric bilinear form, is called an extended affine root supersystem if the root string property is satisfied. $R$ is called a locally finite root supersystem if the form is nondegenerate. Irreducible  locally finite root supersystems have been classified in \cite{you2}.
{One also knows   from \cite{you2} that root string property for a locally finite root supersystem can be replaced by the locally finiteness of the real part.  Generalized root systems according to  the first alternative definition mentioned above, are nothing but finite locally finite root supersystems defined over the complex numbers.} Locally finite root supersystems  naturally appear in the theory of locally finite Lie superalgebras; see  \cite{P} and \cite{you3}.
 Extended affine root supersystems are extensions of locally finite root supersystems by abelian groups and appear as the   root systems of extended affine Lie superalgebras introduced in \cite{you3}; {in particular the root system of an affine Lie superalgebra \cite{van-de} is an extended affine root supersystem.}
The nonzero elements of an extended affine root supersystem are divided into three disjoint parts: One  consists of all real  roots, i.e., the elements  which are not self-orthogonal. The second part is the intersection of the radical of the form with the nonzero elements; the elements of this part are called isotropic  roots. The last part consists of the elements which are not neither isotropic nor real and referred to as nonsingular roots.  An extended affine root supersystem with no nonsingular root is called an affine reflection system  \cite{N1} and an affine reflection system with no isotropic root is called a locally finite root system \cite{LN}.

The concept of a base is so important in the   theory of affine reflection systems and the corresponding Lie algebras. More precisely,   reflectable bases are important in the study of the  structure of locally extended affine root systems \cite{yos1} and integral bases  are important in the theory of  locally finite Lie algebras \cite{NS}.
 A linearly independent subset $\Pi$ of the set of   real  roots of an  affine reflection system  is called  a reflectable base   if  all nonzero reduced  real roots can be obtained from  the iterated  action  of reflections based on the elements of $\Pi.$ Reflectable bases for affine reflection systems have been studied in \cite{AYY}.
A linearly independent  subset  $\Pi$ of a locally finite root supersystem $R$ is called an integral  base if each element of  $R$ can be written as a $\bbbz$-linear combination of the elements of $\Pi.$

{In this work, we  first derive some generic  properties of extended affine root supersystems and locally finite root supersystems and then describe the structure of extended affine root supersystems}.   It is immediate from our results that an irreducible locally finite root supersystem can be recovered from  a nonzero nonsingular root together with a  reflectable base of the real part using  the iterated  action of reflections. We also show that each  locally finite root supersystem $R$  possesses an integral base  and that if $R$ is infinite, then  it has an integral base $\Pi$ with the property that   each element of $R\setminus\{0\}$  can be written as   $r_1\a_1+\cdots+r_n\a_n$ in which $r_1,\ldots,r_n\in\{\pm1\}$ and $\{\a_1,\ldots,\a_n\}\sub \Pi$ with  $r_1\a_1+\cdots+r_t\a_t\in R$ for all $1\leq t\leq n.$
{The result of this paper forms a framework to  study the  locally finite basic classical simple Lie superalgebras \cite{you4}}.

\section{Generic Properties}
Throughout this work, $\bbbf$ is a field of characteristic zero. Unless otherwise mentioned, all vector spaces are considered over $\bbbf.$
We  denote the dual space of a  vector space $V$ by  $V^*.$  We denote the degree of a homogenous element $u$ of a superspace by $|u|$ and make a convention that if  in an expression, we use $|u|$ for an element $u$ of a superspace, by default we have assumed $u$ is homogeneous.    We denote the group of automorphisms of an abelian  group $A$ or a Lie superalgebra $A$   by $Aut(A)$ and for a subset $S$ of an abelian group, by $\la S\ra,$ we mean the subgroup generated by $S.$
 For a set $S,$ by $|S|,$ we mean the cardinal number of $S.$ For a map $f:A\longrightarrow B$ and $C\sub A,$ by $f\mid{_C},$ we mean the restriction of $f$ to $C.$ For two symbols $i,j,$ by $\d_{i,j},$ we mean the Kronecker delta, also $\biguplus$ indicates the disjoint union.
We finally recall that the direct union is, by definition,  the direct limit of a direct system
whose morphisms are inclusion maps.

In the sequel, by a {\it symmetric form} (with values in $\bbbf$) on an additive abelian group $A,$ we mean a map $\fm: A\times A\longrightarrow \bbbf$ satisfying
\begin{itemize}
\item $(a,b)=(b,a)$ for all $a,b\in A,$
\item $(a+b,c)=(a,c)+(b,c)$ and $(a,b+c)=(a,b)+(a,c)$ for all $a,b,c\in A.$
\end{itemize}
In this case, we set  $A^0:=\{a\in A\mid(a,A)=\{0\}\}$ and call it the {\it radical} of the form $\fm.$ The form is called {\it nondegenerate} if $A^0=\{0\}.$
We note that if the form is nondegenerate, $A$ is torsion free and we can identify $A$ as a subset of $\bbbq\ot_\bbbz A.$ In the following, if an abelian group  $A$ is  equipped with a nondegenerate symmetric form, we consider $A$ as a subset of $\bbbq\ot_\bbbz A$ without further explanation.
Also if $V$  is a vector space over a subfield  $\bbbk$ of $\bbbf,$   by a {\it symmetric bilinear form} (with values in $\bbbf$) on $V,$ we mean a map $\fm: V\times V\longrightarrow \bbbf$ satisfying
\begin{itemize}
\item $(a,b)=(b,a);$ \; $(a,b\in V),$
\item $(ra+b,c)=r(a,c)+(b,c)$ and $(a,rb+c)=r(a,b)+(a,c);$\; $(a,b,c\in V,r\in \bbbk).$
\end{itemize}
We set  $V^0:=\{a\in V\mid(a,V)=\{0\}\}$ and call it the {\it radical} of the form $\fm.$ The form is called {\it nondegenerate} if $V^0=\{0\}.$
{We draw the  attention of the  readers to the point that for a $\bbbk$-vector space $V$ equipped with a  symmetric bilinear form $\fm$ with values in $\bbbf$ and a subgroup $A$ of $V,$ the nondegeneracy of $\fm$ and $\fm\mid_{A\times A}$  are not necessarily equivalent.
}

{In this section, we first define an extended affine root supersystem and then study some generic properties of  extended affine root supersystems. Proposition \ref{lfrss} and Lemma \ref{base} are used to get the structure of extended affine root supersystems.
}
\begin{Def}\label{iarr}
{\rm Suppose that $A$ is a nontrivial additive abelian group, $R$ is a subset of $A$ and  $\fm:A\times A\longrightarrow \bbbf$ is  a symmetric form. Set
$$\begin{array}{l}
R^0:=R\cap A^0,\\
\rcross:=R\setminus R^0,\\
\rcross_{re}:=\{\a\in R\mid (\a,\a)\neq0\},\;\;\;\rre:=\rcross_{re}\cup\{0\},\\
\rcross_{ns}:=\{\a\in R\setminus R^0\mid (\a,\a)=0 \},\;\;\; \rim:=\rcross_{ns}\cup\{0\}.
\end{array}$$
We say $(A,\fm,R)$ is an {\it extended affine root supersystem}  if the following hold:
$$\begin{array}{ll}
(S1)& \hbox{$0\in R$ and $\la R \ra= A,$}\\\\
(S2)& \hbox{$R=-R,$}\\\\
(S3)&\hbox{for $\a\in \rre^\times$ and $\b\in R,$ $2(\a,\b)/(\a,\a)\in\bbbz,$}\\\\
(S4)&\parbox{4.5in}{ ({\it root string property}) for $\a\in \rre^\times$ and $\b\in R,$  there are nonnegative  integers  $p,q$  with $2(\b,\a)/(\a,\a)=p-q$ such that \begin{center}$\{\b+k\a\mid k\in\bbbz\}\cap R=\{\b-p\a,\ldots,\b+q\a\};$\end{center}we call  $\{\b-p\a,\ldots,\b+q\a\}$ the {\it $\a$-string through $\b,$}
} \\\\
(S5)&\parbox{4.5in}{for $\a\in \rim$ and $\b\in R$ with $(\a,\b)\neq 0,$
$\{\b-\a,\b+\a\}\cap R\neq \emptyset.$ }
\end{array}$$
If there is no confusion, for the sake of simplicity, we say   {\it $R$ is an extended affine root supersystem in $A.$}
Elements of $R^0$ are called {\it isotropic roots,} elements of  $\rre$ are  called  {\it real roots} and  elements of $\rim$ are  called  {\it nonsingular roots}.
A subset $X$ of $R^\times$ is called {\it connected} if each two elements $\a,\b\in X$ are connected in $X$ in the sense that there is a chain $\a_1,\ldots,\a_n\in X$ with $\a_1=\a,$ $\a_n=\b$ and $(\a_i,\a_{i+1})\neq0,$  $i=1,\ldots,n-1.$ We say an  extended affine root supersystem $R$ is  {\it irreducible} if  $\rre\neq\{0\}$ and $\rcross$ is connected (equivalently, $R^\times$ cannot be written as a disjoint union of two nonempty orthogonal subsets) and say {it is  {\it tame} if for each $\a\in R^0,$ there is $\b\in R^\times$ such that $\a+\b\in R.$}
An extended affine root supersystem $(A,\fm,R)$ is called a {\it locally finite root supersystem} if the form $\fm$ is nondegenerate} and it is called an {\it  affine reflection system} if $\rim=\{0\}.$
\end{Def}

\begin{Examp}
{\rm (1) Suppose that $\LL$ is a finite dimensional basic classical simple Lie superalgebra with a Cartan subalgebra of the even part  and corresponding root system $R.$ One gets from the   finite dimensional  Lie superalgebra theory  that $R$  is a locally finite root supersystem; see \cite{sch}.

(2) {Suppose that $\LL$ is a contragredient Lie superalgebra of finite growth with symmetrizable Cartan matrix \cite{van-de}, then the corresponding root system is an extended affine root supersystem; see \cite[Exa. 3.4 \& Cor. 3.9]{you3}.} }
\end{Examp}

\begin{lem}\label{nonzeero-im}
Suppose that $(A,\fm,R)$ is an extended affine  root supersystem.

(i) If $\a\in \rre$ and $\d\in\rim$ with $(\d,\a)\neq 0,$ then there is a  unique $r\in\{\pm1\}$ such that $\d+r\a\in R.$

(ii) If $\d\in\rim^\times,$ then there is $\eta\in\rim$ with $(\d,\eta)\neq0.$
\end{lem}
\pf
$(i)$ By (S5),  there is  $r\in\{\pm1\}$ such that $\d+r\a\in R.$
  Suppose to the contrary that for $r,s$ with $\{r,s\}=\{1,-1\},$ we have $\b:=\d+s\a,\gamma:=\d+r\a\in R.$ Since $(\b,\d),(\gamma,\d)\neq0,$ we get $\b,\gamma\not\in R^0.$ Also we know that  at most  one of the roots $\b,\gamma$
can be  a nonsingular  root. Suppose that $\b$ is a nonzero  real root, then $(\b,\b)\neq0$ and so $m:=\frac{2s(\d,\a)}{(\a,\a)}\in\bbbz\setminus\{-1\}.$ Since $\b\in \rre^\times,$ we have
\begin{eqnarray*}
  \frac{m}{1+m}=\frac{2s(\d,\a)/(\a,\a)}{1+2s(\d,\a)/(\a,\a)}=\frac{2s(\d,\a)}{(\a,\a)+2s(\d,\a)}&=&\frac{2(\d,\d+s\a)}{(\d+s\a,\d+s\a)}\\
  &=&2(\d,\b)/(\b,\b)\in\bbbz.
\end{eqnarray*}
This implies that $m=-2.$ Now considering the $s\a$-string through $\d,$ we find nonnegative integers $p,q$ with $p-q=-2$ such that $\{\d+ks\a\mid k\in\bbbz \}\cap R=\{\d-ps\a,\ldots,\d+qs\a\};$ in particular as $\d-s\a=\d+r\a=\gamma\in R,$ we have $\d+3s\a\in R.$ But $$(\d+3s\a,\d+3s\a)=6s(\d,\a)+9(\a,\a)=-6(\a,\a)+9(\a,\a)=3(\a,\a)\neq0$$ and $$\frac{2(\a,\d+3s\a)}{(\d+3s\a,\d+3s\a)}=\frac{2(\a,\d)+6s(\a,\a)}{3(\a,\a)}=\frac{-2s(\a,\a)+6s(\a,\a)}{3(\a,\a)}
=\frac{4s}{3}\not\in\bbbz,$$ a contradiction. This completes the proof.

$(ii)$ Since $\d\in\rim^\times,$ we have   $\d\not\in A^0.$ Therefore,  there is $\a\in R^\times$ with $(\d,\a)\neq0.$ If $\a$ is nonsingular, we are done, so suppose $\a\in \rre^\times.$ Set $n:=\frac{2(\d,\a)}{(\a,\a)}\in\bbbz.$  Considering the $\a$-string through $\d,$ we find nonnegative integers  $p,q$ with $p-q=n$ such that $\{k\in\bbbz\mid \d+k\a\in R\}=\{-p,\ldots,q\}.$ Since $-p\leq -n\leq q,$  we have $\eta:=\d-n\a\in R.$ Now we have $(\d,\eta)=(\d,\d-n\a)=-n(\d,\a)\neq0$ and $(\eta,\eta)=(\d-n\a,\d-n\a)=n^2(\a,\a)-2n(\d,\a)=0.$ So $\eta\in \rim $ with $(\d,\eta)\neq0.$
\qed

\begin{lem}\label{abr}
Suppose that $A$ is a nontrivial additive abelian group, $R$ is a subset of $A$ and  $\fm:A\times A\longrightarrow \bbbf$ is  a nondegenerate symmetric form. If  $(A,\fm,R)$ satisfies $(S1),(S3)-(S5),$ then $(S2)$ is also satisfied.
\end{lem}
\pf  We assume $\a\in R.$ We must prove that $-\a\in R.$ If $\a\in \rre^\times,$ then the root string property implies that $\a-2\a\in R$ and so $-\a\in R.$ Next suppose that $\a\in \rim^\times,$ then using the same argument as in Lemma \ref{nonzeero-im}($ii$), we find $\eta\in \rim$ with $(\a,\eta)\neq0.$
So there is $r\in\{\pm 1\}$ with $\b:=\a+r\eta\in R.$ Since $\b\in \rre,$ we have $-\b\in \rre.$  On the other hand, $(-\b,\eta)\neq0,$ so either $-\b+r\eta\in R$ or $-\b-r\eta\in R.$  But if  $-\b-r\eta=-\a-2r\eta\in R,$ we get $-\a-2r\eta\in \rre^\times$ while $$2\frac{(\eta,-\b-r\eta)}{(-\b-r\eta,-\b-r\eta)}=2\frac{(\eta,-\a-2r\eta)}{(-\a-2r\eta,-\a-2r\eta)}=-r/2\not\in\bbbz$$ which is a contradiction. So $-\a=-\b+r\eta\in R.$\qed
\begin{Def}
{\rm Suppose that $(A,\fm, R)$ is a locally finite root supersystem.
\begin{itemize}
\item
 The subgroup $\w$ of $Aut(A)$ generated by $r_\a$ ($\a\in\rre^\times$) mapping  $a\in A$ to $a-\frac{2(a,\a)}{(\a,\a)}\a,$  is called the {\it Weyl group} of $R.$
\item A subset $S$ of  $R$ is called a  {\it sub-supersystem} if the restriction of the form to $\la S\ra$ is nondegenerate, $0\in S,$ for $\a\in S\cap\rre^\times, \b\in S$ and $\gamma\in S\cap\rim$ with $(\b,\gamma)\neq 0,$ $r_\a(\b)\in S$ and  $\{\gamma-\b,\gamma+\b\}\cap
    S\neq\emptyset.$
\item A sub-supersystem $S$ of $R$ is called  {\it $\bbbz$-linearly  closed} if $R\cap (\hbox{span}_\bbbz S)=S.$
\item If  $(A,\fm, R)$ is irreducible, $R$ is said to be  of {\it real type} if  $\hbox{span}_\bbbq R_{re}=\bbbq\ot_{\bbbz} A;$ otherwise, we say it is of {\it imaginary type.}
\item  If $\{R_i\mid i\in I\}$ is a class of  sub-supersystems of $R$ which are mutually   orthogonal with respect the form $\fm$ and $R\setminus\{0\}=\uplus_{i\in I}(R_i\setminus\{0\}),$  we say $R$ is {\it the direct sum} of $R_i$'s and write $R=\op_{i\in I}R_i.$
\item The locally finite  root supersystem $(A,\fm,R)$ is called a {\it locally finite root system} if $\rim=\{0\}.$
\item $(A,\fm, R)$ is  said to be {\it isomorphic} to another  locally finite root supersystem $(B,\fm',S)$ if there is a group isomorphism $\varphi:A\longrightarrow B$ and a nonzero scalar $r\in\bbbf$ such that $\varphi(R)=S$ and  $(
a_1,a_2)=r(\varphi(a_1),\varphi(a_2))'$ for all $a_1,a_2\in A.$
\end{itemize}}
\end{Def}

\begin{Rem}\label{rem2}
{\rm
(i)  Locally finite root systems  initially  appeared in the work of K.H. Neeb and N. Stumme \cite{NS} on locally finite split simple Lie algebras. Then in 2003, O. Loos and E. Neher \cite{LN} systematically studied locally finite root systems. In their sense a locally finite root system is a locally finite spanning set $R$ of a nontrivial  vector space $\v$ such that $0\in R$ and for each  $\a\in R\setminus\{0\},$ there is a functional $\check\a\in\v^*$ such that $\check\a(\a)=2,$ $\check\a(\b)\in\bbbz$ for all $\b\in R$ and that $\b-\check\a(\b)\a\in R.$ It is proved that locally finiteness can be replaced by the existence of a nonzero bilinear form which is positive definite on the $\bbbq$-span of $R$ and invariant under the Weyl group; moreover such a form is nondegenerate and is unique up to a  scalar multiple if $R$ is irreducible \cite[\S 4.1]{LN}. Also a locally finite root system  $R$ in $\v$ contains a $\bbbz$-basis for $\la R\ra$ \cite[Lem. 5.1]{LN2}. This allows us to have a natural isomorphism between $\v$ and $\bbbf\ot_\bbbz\la R\ra$ and so it is natural to consider a locally finite root system as a subset of a torsion free abelian group instead of a subset of a vector space.

(ii) Suppose that $S$ is a  sub-supersystem of a locally finite root supersystem $R,$ then $S_{re}$ is a locally finite root system by \cite[\S 3.1]{you2} and \cite[\S 3.4]{LN}. Now the same argument as in \cite[Lem. 3.12]{you2} shows that the root string property holds for $S.$ This together with Lemma \ref{abr} implies that $S$ is a locally finite root supersystem in its $\bbbz$-span.
}
\end{Rem}

Suppose that $T$ is a nonempty  index set  with $|T|\geq 2$ and $\u:=\op_{i\in
T}\bbbz\ep_i$ is the free $\bbbz$-module over   the
set $T.$ Define the  form $$\begin{array}{c}\fm:\u\times\u\longrightarrow\bbbf\\
(\ep_i,\ep_j)\mapsto\d_{i,j}, \hbox{ for } i,j\in T,
\end{array}$$
and set
\begin{equation}\label{locally-finite}
\begin{array}{l}
\dot A_T:=\{\ep_i-\ep_j\mid i,j\in T\},\\
D_T:=\dot A_T\cup\{\pm(\ep_i+\ep_j)\mid i,j\in T,\;i\neq j\},\\
B_T:=D_T\cup\{\pm\ep_i\mid i\in T\},\\
C_T:=D_T\cup\{\pm2\ep_i\mid i\in T\},\\
BC_T:=B_T\cup C_T.
\end{array}
\end{equation}
These are irreducible locally finite root systems
in their $\bbbz$-span's. Moreover, each irreducible  locally finite root system is either  an irreducible finite root system  or a locally finite root system  isomorphic to one of these locally finite root
systems. We refer to locally finite root systems listed in (\ref{locally-finite}) as  {\it type} $A,D,B,C$
and $BC$  respectively. We note that if $R$ is  an irreducible locally finite
root system as above, then  $(\a,\a)\in\bbbn$ for all
$\a\in R.$ This allows us to  define
$$\begin{array}{l}
R_{sh}:=\{\a\in R^\times\mid (\a,\a)\leq(\b,\b);\;\;\hbox{for all $\b\in R$} \},\\
R_{ex}:=R\cap2 R_{sh}\andd
R_{lg}:= R^\times\setminus( R_{sh}\cup R_{ex}).
\end{array}$$
The elements of $R_{sh}$ (resp. $R_{lg},R_{ex}$) are called {\it
short roots} (resp. {\it long roots, extra-long roots}) of $R$. We point  out that following the usual notation in the literature,  the locally finite root system  of type $A$  is denoted by  $\dot A$ instead of $A,$ as all locally finite root systems listed above are spanning sets for $\bbbf\ot_\bbbz \u$ other than the one of type $A$ which spans a subspace of codimension 1.

\begin{lem}\label{super-sys}

(i) If $\{(X_i,\fm_i,S_i)\mid i\in I\}$ is a class of   locally finite root supersystems, then for   $X:=\op_{i\in I}X_i$ and $\fm:=\op_{i\in I}\fm_i,$   $(X,\fm,S:=\cup_{i\in I} S_i)$ is a locally finite root supersystem.

(ii) Suppose that $( A,\fm, R)$ is a locally finite root supersystem.  Connectedness is an equivalence relation on $ R\setminus\{0\}.$ Also if $S$ is a connected component of $ R\setminus\{0\},$ then $S\cup\{0\}$ is an irreducible sub-supersystem of $ R.$ Moreover, $ R$ is a direct sum of irreducible sub-supersystems.

(iii) Suppose that $( A,\fm, R)$ is a locally finite root supersystem. For $ A_{re}:=\la R_{re}\ra$ and $\fm_{re}:=\fm\mid_{_{ A_{re}\times  A_{re}}},$ $( A_{re},\fm_{re}, R_{re})$ is a locally finite root system.
%
\end{lem}
\pf  See \cite[\S 3]{you2}.
\qed

\medskip

{We also have the following straightforward lemma:}
\begin{lem} \label{direct-sum}
Suppose that $(A,\fm,R)$ is  an irreducible  locally finite root supersystem, set  $\v:=\bbbf\ot_\bbbz A$ and identify $A$ as a subset of $\v.$  Then  $\v=\hbox{span}_\bbbf \rre$ if and only if $R$ is of real type.
\end{lem}

In the following two theorems, we give the classification of irreducible locally finite root supersystems.

\begin{The}[{\cite[Thm. 4.28]{you2}}]\label{classification I}
Suppose that $T,T'$ are index sets of cardinal numbers greater than $1$ with  $|T|\neq|T'|$ if $T,T'$ are both finite. Fix a symbol $\a^*$ and  pick $t_0\in T$ and $p_0\in T'.$ Consider the free $\bbbz$-module $X:=\bbbz\a^*\op\op_{t\in T}\bbbz\ep_t\op\op_{p\in T'}\bbbz\d_p$  and define the symmetric form $$\fm:X\times X\longrightarrow \bbbf$$ by $$\begin{array}{ll}
(\a^*,\a^*):=0,(\a^*,\ep_{t_0}):=1,(\a^*,\d_{p_0}):=1\\
(\a^*,\ep_t):=0,(\a^*,\d_q):=0&t\in T\setminus\{t_0\},q\in T'\setminus\{p_0\}\\
(\ep_t,\d_p):=0,(\ep_t,\ep_s):=\d_{t,s},(\d_p,\d_q):=-\d_{p,q}&t,s\in T,p,q\in T'.
\end{array}$$
 Take $R$ to be $\rre\cup \rim^\times$ as in the following table:
$${\small
\begin{tabular}{|c|c|c|}
\hline
type &$\rre$&$\rim^\times$
\\\hline
$\dot A(0,T)$& $\{\ep_t-\ep_s\mid t,s\in T\}$&$\pm\w\a^*$\\
\hline
$\dot C(0,T)$& $\{\pm(\ep_t\pm\ep_s)\mid t,s\in T\}$&$\pm\w\a^*$\\
\hline
$\dot A(T,T')$& $\{\ep_t-\ep_s,\d_p-\d_q\mid t,s\in T,p,q\in T'\}$&$\pm\w\a^*$\\
\hline
\end{tabular}}$$
in which  $\w$ is the subgroup of $Aut(X)$ generated by the reflections $r_\a$ $(\a\in \rre\setminus\{0\})$ mapping $\b\in X$ to $\b-\frac{2(\b,\a)}{(\a,\a)}\a,$ then $(A:=\la R\ra,\fm\mid_{A\times A}, R)$  is an irreducible locally finite root supersystem of imaginary type and conversely, each irreducible locally finite root supersystem of imaginary type  is isomorphic to one and only one of these root supersystems.
\end{The}

\begin{The}[{\cite[Thm. 4.37]{you2}}]
\label{classification II}
 Suppose $n\in \{2,3\}$ and  $(X_1,\fm_1,S_1),$ $\ldots,(X_n,\fm_n,S_n),$  are irreducible locally finite root systems. Set $X:=X_1\op\cdots\op X_n$ and $\fm:=\fm_1\op\cdots\op\fm_n$ and consider  the locally finite root system  $(X,\fm,S:=S_1\op \cdots\op S_n).$ Take $\w$ to be the Weyl group of $S.$  For $1\leq i\leq n,$ we identify $X_i$ with a subset of $\bbbq\ot_\bbbz X_i$ in the usual manner. If $1\leq i\leq n$ and $S_i$ is a finite root system of rank $\ell\geq 2,$ we take $\{\omega_1^i,\ldots,\omega_\ell^i\}\sub\bbbq\ot_\bbbz X_i$ to be a  set of fundamental weights for $S_i$ (see \cite[Pro. 2.7]{you2}) and if $S_i$ is one of infinite locally finite root systems $B_T, C_T, D_T$ or $BC_T$  as in (\ref{locally-finite}), by $\omega_1^i,$ we mean $\ep_1,$ where $1$ is a distinguished element of $T.$ Also if $S_i$ is  one of the  finite root systems $\{0,\pm\a\}$ of type $A_1$ or  $\{0,\pm\a,\pm2\a\}$ of type $BC_1,$ we set $\omega_1^i:=\frac{1}{2}\a.$
Consider $\d^*$ and $ R:= R_{re}\cup R_{ns}^\times$ as in the following table:
$${\tiny
\begin{tabular}{|c|l|c|c|c|c|}
\hline
$ n$& $ S_i\;(1\leq i\leq n)$&$ R_{re}$&$\d^*$&$ R_{ns}^\times$&\hbox{type}\\
\hline
$ 2$& $ S_1=A_\ell,\; S_2=A_\ell$ $(\ell\in\bbbz^{\geq1})$&$S_1\op S_2$&$\omega_1^1+\omega_1^2$&$\pm\w\d^*$&$A(\ell,\ell)$\\
\hline
$2$& $ S_1=B_T,\;S_2=BC_{T'}$ $(|T|,|T'|\geq2)$&$S_1\op S_2$&$\omega_1^1+\omega_1^2$&$\w\d^*$&$B(T,T')$\\
\hline
$2$& $ S_1=BC_T,\;S_2=BC_{T'}$ $(|T|,|T'|>1)$&$S_1\op S_2$&$\omega_1^1+\omega_1^2$&$\w\d^*$&$BC(T,T')$\\
\hline
$2$& $ S_1=BC_T,\;S_2=BC_{T'}$ $(|T|=1,|T'|=1)$&$S_1\op S_2$&$2\omega_1^1+2\omega_1^2$&$\w\d^*$&$BC(T,T')$\\
\hline
$2$& $ S_1=BC_T,\;S_2=BC_{T'}$ $(|T|=1,|T'|>1)$&$S_1\op S_2$&$2\omega_1^1+\omega_1^2$&$\w\d^*$&$BC(T,T')$\\
\hline
$2$& $ S_1=D_T,\;S_2=C_{T'}$ $(|T|\geq3,|T'|\geq2)$&$S_1\op S_2$&$\omega_1^1+\omega_1^2$&$\w\d^*$&$D(T,T')$\\
\hline
$2$& $ S_1=C_T,\;S_2=C_{T'}$ $(|T|,|T'|\geq2)$&$S_1\op S_2$&$\omega_1^1+\omega_1^2$&$\w\d^*$&$C(T,T')$\\
\hline
$2$& $ S_1=A_1,\;S_2=BC_{T}$ $(|T|=1)$&$S_1\op S_2$&$2\omega_1^1+2\omega_1^2$&$\w\d^*$&$B(1,T)$\\
\hline
$2$& $ S_1=A_1,\;S_2=BC_{T}$ $(|T|\geq2)$&$S_1\op S_2$&$2\omega_1^1+\omega_1^2$&$\w\d^*$&$B(1,T)$\\
\hline
$2$& $ S_1=A_1,\;S_2=C_T$ $(|T|\geq 2)$&$S_1\op S_2$&$\omega_1^1+\omega_1^2$&$\w\d^*$&$C(1,T)$\\
\hline
$2$& $ S_1=A_1,\;S_2=B_3$ &$S_1\op S_2$&$\omega_1^1+\omega_3^2$&$\w\d^*$&$AB(1,3)$\\
\hline
$2$& $ S_1=A_1,\;S_2=D_{T}$ $(|T|\geq3)$&$S_1\op S_2$&$\omega_1^1+\omega_1^2$&$\w\d^*$&$D(1,T)$\\
\hline
$2$& $ S_1=BC_1,\;S_2=B_T$ $(|T|\geq2)$&$S_1\op S_2$&$2\omega_1^1+\omega_1^2$&$\w\d^*$&$B(T,1)$\\
\hline
$2$& $ S_1=BC_1,\;S_2=G_2$ &$S_1\op S_2$&$2\omega_1^1+\omega_1^2$&$\w\d^*$&$G(1,2)$\\
\hline
$3$& $ S_1=A_1,\;S_2=A_1,\; S_3=A_1$ &$S_1\op S_2\op S_3$&$\omega_1^1+\omega_1^2+\omega_1^3$&$\w\d^*$&$D(2,1,\lam) {\tiny (\lam\neq 0,-1)}$\\
\hline
$3$& $ S_1=A_1,\; S_2= A_1,\; S_3:=C_T$ $(|T|\geq2)$&$S_1\op S_2\op S_3$&$\omega_1^1+\omega_1^2+\omega_1^3$&$\w\d^*$&$D(2,T)$\\
\hline
\end{tabular}
}$$
For $1\leq i\leq n,$ normalize the form $\fm_i$ on $X_i$ such that  $(\d^*,\d^*)=0$ and that for type $D(2,T),$ $(\omega_1^1,\omega_1^1)_1=(\omega_1^2,\omega_1^2)_2.$  Then $(\la  R\ra,\fm\mid_{\la R\ra\times \la R\ra}, R)$ is an irreducible locally finite root supersystem of real type and conversely, if  $( X,\fm, R)$ is an irreducible locally finite root supersystem of real type, {it is either an irreducible locally finite root system or  isomorphic to one  of the locally finite root supersystems listed in the above table. Moreover, locally finite root supersystems in the above table are mutually non-isomorphic except for the ones of type $D(2,1,\lam).$ More precisely,  For $\lam,\mu\in\bbbf\setminus\{0,-1\},$ $D(2,1,\lam)$ is isomorphic to $D(2,1,\mu)$  if and only if $\lam,\mu$ are in the same orbit under the action of the group of permutations on $\bbbf\setminus\{0,-1\}$ generated by $\a\mapsto \a^{-1}$ and $\a\mapsto -1-\a.$}
\end{The}

We make  a convention  that from now on for the types listed in column ``type'' of Theorems \ref{classification I} and \ref{classification II}, we may use  a finite  index set $T$ and its cardinal number in place of each other, e.g.,  if $T$ is  a nonempty finite set of cardinal number $\ell,$ instead of type $B(1,T),$ we may write $B(1,\ell).$

\begin{Pro}\label{lfrss}
Suppose that $(A,\fm, R)$ is an extended affine root supersystem and $\bar{\;}:A\longrightarrow \bar A:=A/A^0$  is the {canonical  epimorphism}. Suppose that $\fm^{\bar{\;}}$ is the induced form on $\bar A$ defined by $$(\bar a,\bar b):=(a,b);\;\;\; (a,b\in A).$$ Then we have the following:

$(i)$  $\{2(\a,\b)/(\a,\a)\mid \a\in \rre^\times,\b\in R\}$ is a bounded subset of $\bbbz$ and  for $\a\in R_{re}^\times$ and $\b\in R_{ns},$ $2(\a,\b)/(\a,\a)\in\{0,\pm1,\pm2\}.$


$(ii)$ If $\a,\b\in\rre^\times$ are connected in $\rre,$ then $(\a,\a)/(\b,\b)\in\bbbq.$
Also each subset of $\rre^\times$ whose elements are mutually disconnected in $\rre^\times$ is $\bbbz$-linearly independent.

$(iii)$ $(\bar A,\fm^{\bar{\;}},\bar R)$ is a locally finite root supersystem.  Moreover, if $R$ is irreducible, then so is  $\bar R.$
\end{Pro}
\pf  $(i)$ See \cite[Lem. 3.7]{you2} and follow the proof of \cite[Lem.  3.8]{you2}.


$(ii)$ See \cite[Lem. 3.6]{you2}.

$(iii)$
Set $\v:=\bbbf\ot_{_\bbbz}\bar A.$ Since $\bar A$ is torsion free, we identify $\bar A$ as a subset of $\v$ and set  $\v_\bbbq:=\hbox{span}_\bbbq \bar R$ as well as      $\v_{re}:=\hbox{span}_\bbbq\bar R_{re}.$  The nondegenerate form $\fm^{\bar\;}:\bar A\times \bar A\longrightarrow \bbbf$ induces a bilinear form $$\begin{array}{c}\fm_\bbbf:(\bbbf\ot_\bbbz\bar A)\times(\bbbf\ot_\bbbz\bar A)\longrightarrow \bbbf\\
(r\ot \bar a,s\ot \bar b):=rs(a,b);\;\;(r,s\in\bbbf,\;a,b\in A).\end{array}$$ Take $\fm_\bbbq$ to be the  restriction of the  form $\fm_{\bbbf}$ to $\v_\bbbq=\hbox{span}_\bbbq\bar R.$   Using the same argument as in \cite[Lem. 1.6]{AYY}, one can see that $\fm_\bbbq$ is nondegenerate. To carry out  the proof, we just need to verify the root string property. To this end using \cite[Lem.'s  3.10 \& 3.12]{you2}, it is enough to show that  $ \bar R_{re}= \overline\rre=\{
\bar\a\mid \a\in \rre\}\sub\bbbf\ot\bar A$ is locally finite in $\v_{re}=\hbox{span}_\bbbq \bar R_{re}$ in the sense that it intersects each finite dimensional subspace of $\v_{re}$ in a finite set.
Now we assume $\w$ is a finite dimensional subspace of  $\v_{re}$ and show that $\bar R_{re}\cap \w$ is a finite set. Since $\w$ is a finite dimensional subspace of $\v_{re},$ there is a finite subset $\{\a_1,\ldots,\a_m\}\sub  R_{re}$ such that $\w\sub U_1:=\hbox{span}_\bbbq\{ \bar\a_1,\ldots,  \bar\a_m\}.$ By  \cite[Lem. 3.1]{you2}, there is a finite dimensional subspace $U_2$ of $\v_\bbbq$ such that $U_1\sub U_2$ and  the form $\fm_\bbbq$ restricted to $U_2$ is nondegenerate.
Suppose that $\{R_i\mid i\in I\}$ is the class of connected components of $\rre^\times.$ To complete the proof  {using part ($ii$)} together with the fact that $U_2$ is finite dimensional,
we need   to  show that for all $i\in I,$  $U_2\cap   \bar R_i$ is a finite set.
Since $U_2$ is finite dimensional, there is a finite set $\{\b_1,\ldots,\b_n\}\sub R$ such that $U_2\sub\hbox{span}_\bbbq\{ \bar \b_1,\ldots, \bar \b_n\}.$
Fix $i\in I$ and consider the map $$\begin{array}{l}\varphi:U_2\cap   \bar R_i\longrightarrow \bbbz^n\\
\bar\a\mapsto (\frac{2(\bar\a,\bar\b_1)}{(\bar\a,\bar\a)},\ldots,\frac{2(\bar\a,\bar\b_n)}{(\bar\a,\bar\a)}).\end{array}$$
We claim that $\varphi $ is one to one. Suppose that for  $\a,\b\in R_i,$ $ \bar\a, \bar\b\in U_2\cap   \bar R_i$  and $$(\frac{2(\bar\a,\bar\b_1)}{(\bar\a,\bar\a)},\ldots,\frac{2(\bar\a,\bar\b_n)}{(\bar\a,\bar\a)})=(\frac{2(\bar\b,\bar\b_1)}{(\bar\b,\bar\b)},\ldots,\frac{2(\bar\b,\bar\b_n)}{(\bar\b,\bar\b)}),$$ then for $1\leq i\leq n,$ $\frac{(\bar\a,\bar\b_i)}{(\bar\a,\bar\a)}=\frac{(\bar\b,\bar\b_i)}{(\bar\b,\bar\b)}.$ So $(\frac{  \bar\a}{(\bar\a,\bar\a)}-\frac{  \bar\b}{(\bar\b,\bar\b)}, \bar\b_i)_\bbbf=0$ for all $1\leq i\leq n.$ Therefore,
$(\frac{  \bar\a}{(\bar\a,\bar\a)}-\frac{  \bar \b}{(\bar\b,\bar\b)},U_2)_\bbbf=\{0\}.$  But $\frac{(\bar\a,\bar\a)}{(\bar\b,\bar\b)}\in\bbbq$ (see part $(ii)$) and so $$(  \bar\a-\frac{(\bar\a,\bar\a)}{(\bar\b,\bar\b)}  \bar\b,U_2)_\bbbq=(  \bar\a-\frac{(\bar\a,\bar\a)}{(\bar\b,\bar\b)}  \bar\b,U_2)_\bbbf=\{0\}.$$
So  we get that $\bar\a=\frac{(\bar\a,\bar\a)}{(\bar\b,\bar\b)}\bar\b$ as the form $\fm_\bbbq$ on $U_2$ is nondegenerate.
But as $\frac{2(\bar\a,\bar\b)}{(\bar\a,\bar\a)},\frac{2(\bar\a,\bar\b)}{(\bar\b,\bar\b)}\in \bbbz,$ we get that $(\bar\a,\bar\a)/(\bar\b,\bar\b)\in\{\pm1,\pm 2,\pm\frac{1}{2}\}.$ If $\frac{(\bar\a,\bar\a)}{(\bar\b,\bar\b)}=\pm2,$ then $\bar\a=\pm2\bar\b$ and so $\frac{(\bar\a,\bar\a)}{(\bar\b,\bar\b)}=4,$ a contradiction, also if $\frac{(\bar\a,\bar\a)}{(\bar\b,\bar\b)}=\pm(1/2),$ then $\bar\a=\pm(1/2)\bar\b$ and $(\bar\a,\bar\a)/(\bar\b,\bar\b)=1/4$ which is again  a contradiction. If $(\bar\a,\bar\a)=-(\bar\b,\bar\b),$ then $\bar\a=-\bar\b$ and so $(\bar\a,\bar\a)/(\bar\b,\bar\b)=1$ that is absurd. Therefore,  $\bar\a=\bar\b$ i.e.,  $\varphi$ is one to one. Also  {using part ($i$)}, we get that the set  $\{\frac{2(\bar\a,\bar\b)}{(\bar\a,\bar\a)}\mid \a\in\rre,\;\b\in R\}$ is bounded. This in turn  implies that the image of $\varphi$ and so  $ U_2\cap \bar R_{i}$ is finite. This together with Lemma \ref{super-sys} completes the proof of the first assertion. The last assertion follows from  an immediate verification.  \qed

\begin{Def}{\rm
Suppose that $(A,\fm , R)$ is an irreducible extended affine root supersystem.  We define the \emph{type} of $R$ to be the type of $\bar R.$}
\end{Def}

\begin{lem}\label{radical}
Suppose that $A$ is a torsion free abelian group and $(A,\fm, R)$ is an irreducible extended affine root supersystem of type {\small $X\neq A(\ell,\ell),BC(1,1).$}  Then for each $a\in A^0,$ there is a nonzero integer $n$ such that  $na\in \la R^0\ra;$ in particular, if $ X\neq A(\ell,\ell),$ $R^0=\{0\}$  if and only if $A^0=\{0\}.$
\end{lem}
\pf
Set $\v:=\bbbq\ot_\bbbz A.$ Since $A$ is torsion free, we identify $A$ as  a subset of $\v.$ The form $\fm$ induces the  symmetric bilinear form   $\v\times\v\longrightarrow \bbbf$ (with values in $\bbbf$) defined by  $(r\ot a,s\ot b):=rs(a,b)$ ($r,s\in \bbbq,$ $a,b\in A$); we  denote this bilinear form   again by $\fm.$ Set $\v^0:=\{\a\in\v\mid(\a,\v)=\{0\}\}.$ Suppose that  $\bar{\;}:\v\longrightarrow \bar \v:=\v/\v^0$  is the {canonical epimorphism} and that $\fm^{\bar{\;}}$ is the induced map on $\bar \v\times\bar \v.$
We note that $\v^0=\hbox{span}_\bbbq A^0$ and use Proposition \ref{lfrss} and Lemma \ref{super-sys} to get that   $\bar R_{re}$ is a locally finite root system in its $\bbbz$-span. Therefore by \cite[Lem. 5.1]{LN2}, there is a $\bbbz$-basis $B\sub \bar R_{re}$ for $\bar A_{re}:=\la \bar R_{re}\ra$ such that \begin{equation}\label{ref}\w_BB=(\bar R_{re})^\times_{red}:=\bar R_{re}\setminus\{2\bar \a\mid\a\in R_{re}\},\end{equation} in which by $\w_B,$ we mean  the subgroup of the Weyl group of $\bar R_{re}$ generated by $r_{\bar\a}$ for all $\bar\a\in B.$ Fix $\a^*\in \bar R_{ns}^\times$ if $\bar R$ is of imaginary type and   set $$K:=\left\{\begin{array}{cl}B&\hbox{if $\bar R$ is of real type,}\\
B\cup\{\a^*\}& \hbox{if $\bar R$ is of imaginary type}.
\end{array}\right.$$
Then $K$ is a basis for  $\bbbq$-vector space   $\bar \v.$  Take $\dot K\sub R$ to be a preimage of $K$ under the canonical map `` $\bar{}$ '', then $\dot K$ is a $\bbbq$-linearly independent subspace  of $\v$ and for $\dot \v:=\hbox{span}_\bbbq \dot K,$ we have $\v=\dot\v\op\v^0.$
Now set  $\dot R:=\{\dot\a\in\dot\v\mid \exists \sg\in\v^0; \dot\a+\sg\in R\}$ and for each $\dot \a\in \dot R,$ set $T_{\dot \a}:=\{\sg\in \v^0\mid \dot\a+\sg\in R\}.$ Then $\dot R$ is a locally finite root supersystem in its $\bbbz$-span isomorphic to $\bar R.$ Since $\dot K\sub R\cap\dot R,$ we have $-\dot K\sub R\cap\dot R.$ Taking  $\w_{\dot K}$ to be  the subgroup of the Weyl group of $R$  generated by the reflections based on  real roots of $\dot K,$ we have
$$\w_{\dot K}(\pm {\dot K})\sub R\cap \dot R\andd \pm\w_{\dot K}{\dot K}=\left\{\begin{array}{ll}(\dot R_{re})_{red}^\times&\hbox{if $\bar R$ is of real type,}\\
\dot R^\times& \hbox{if $\bar R$ is of imaginary  type.}\end{array}\right.$$
So \begin{equation}\label{0}\left\{\begin{array}{ll}
0\in T_{\dot\a}& \hbox{if $\dot R$ is of real type and $\dot\a\in (\dot R_{re})_{red}:=(\dot R_{re})_{red}\cup\{0\},$}\\
0\in T_{\dot\a}& \hbox{if $\dot R$ is of imaginary type and $\dot\a\in \dot R.$}
\end{array}\right.\end{equation}

To proceed with  the proof, we claim that    for each $\dot\a\in \dot R$ and $\sg\in T_{\dot\a},$ there is $n\in\bbbz\setminus\{0\}$ such that $n\sg\in \la R^0\ra.$ If $\dot \a=0, T_{\dot\a}\sub R^0$ and there is nothing to prove. Now   the following cases can happen:
\smallskip

\underline{Case 1. $\dot\a\in{\dot R_{re}}^\times:$}  In this case, we show that $T_{\dot\a}\sub R^0.$  We first assume $\dot\a\in (\dot R_{re})_{red}^\times,$ then since $0\in T_{\dot\a},$ $\a:=\dot\a,\b:=\dot\a+\sg\in R.$ Now considering the $\a$-string through $\b,$ we find that $\sg\in R$ and so it is an element of $R^0.$ Next suppose that  $\dot\a\in \dot R_{re}^\times\setminus (\dot R_{re})_{red},$  then  there exists $\dot\b\in (\dot R_{re})_{red}$ with $\dot\a=2\dot\b.$ Now for  $\sg\in T_{\dot\a},$ takeing $\a:=\dot\b$ and $\b:=\dot\a+\sg$ and considering the $\a$-string through $\b,$ we get  that $\sg\in R^0.$
\smallskip

\underline{Case 2. $\dot R$ is of real type and  $\dot\a\in \dot R_{ns}^\times:$} For $\dot \gamma\in (\dot R_{re})_{red}^\times$ and $\eta\in  T_{\dot\a},$ since $\dot\gamma\in \rre^\times,$ we have $r_{\dot\gamma}(\dot\a+\eta)=r_{\dot\gamma}(\dot\a)+\eta\in R.$ This implies that $T_{\dot\a}\sub T_{r_{\dot\gamma}(\dot\a)};$ similarly we have  $T_{r_{\dot\gamma}(\dot\a)}\sub T_{\dot\a}.$ We know that the Weyl group $\dot\w$ of $\dot R$ is generated by the reflections based on nonzero elements of $(\dot R_{re})_{red}$  and that each two nonzero  nonsingular roots are $\dot \w$-conjugate as $\dot R$ is not of type $A(\ell,\ell).$ These altogether imply that $T:=T_{\dot\a}=T_{\dot\b}$ for all nonzero nonsingular  roots $\dot\b.$  Since $\dot R$ is of real type $X\neq  BC(1,1),A(\ell,\ell),$ one finds nonsingular roots $\dot\b,\dot\gamma$ with $(\dot\gamma,\dot\b)\neq0,$ $\dot\b-\dot\gamma\in \dot R_{re}$ and $\dot\b+\dot\gamma\not\in \dot R.$ We next note that $T=T_{\dot\a}=-T_{-\dot\a}=-T$ and fix $\sg,\tau\in T=-T.$ Since $\a:=\dot\b+\sg,\b:=\dot\gamma+\tau,\gamma:=\dot\gamma-\tau\in R$ and $(\a,\b),(\a,\gamma)\neq0,$ there are $r,s\in\{\pm1\}$ with $\zeta:=\a+r\b,\eta:=\a+s\gamma\in  R.$ But $\dot\b+\dot\gamma\not\in \dot R,$  so $\zeta=\dot\b-\dot\gamma+\sg-\tau,\eta=\dot\b-\dot\gamma+  \sg+\tau.$ Therefore using the  previous case, we have  $\sg-\tau,\sg+\tau\in R^0;$ this in particular  implies that $2\sg,2\tau\in \la R^0\ra.$

\smallskip
\underline{Case 3. $\dot R$ is of imaginary type and  $\dot\a\in \dot R_{ns}^\times:$} By \cite[Lem. 4.5]{you2}, there is $\dot\b\in \dot R_{re}$ such that $(\dot\a,\dot\b)\neq 0.$ We next note that     $T:=T_{\dot\a}=-T_{-\dot\a}.$  Also as  $0\in T_{\dot\b}$ and $R$ is invariant under the reflections,   $T:=T_{\dot\a}=T_{r_{\dot\b}(\dot\a)}$ as in the previous case. Also by \cite[Lem.'s 4.6 \&  4.7]{you2}, we have $r_{\dot\b}\dot\a-\dot\a\in \dot R_{re}$  while $r_{\dot\b}\dot\a+\dot\a\not\in \dot R.$   Now for $\sg,\tau\in T,$ we have $(r_{\dot\b}\dot\a+\sg,\dot\a+\tau)\neq 0.$ Since $r_{\dot\b}\dot\a+\dot\a\not\in \dot R,$ we get that $r_{\dot\b}\dot\a-\dot\a+\sg-\tau\in \rre$ and so using Case 1,    we have $\sg-\tau\in R^0.$ Thus we have  $T-T\sub R^0;$ but   $0\in T,$ so  $T=T_{\dot\a}\sub R^0.$

Now suppose $a\in A^0\setminus\{0\},$ then $a\in \v^0$ and there are $r_1,\ldots,r_m\in \bbbz\setminus\{0\}$ and $\a_1,\ldots,\a_m\in R\setminus\{0\}$  with $a=\sum_{i=1}^m r_i\a_i.$ But for each $1\leq i\leq m,$ there are $\dot\a_i\in \dot R,$  $n_i\in \bbbz\setminus\{0\}$ and $\d_i\in \la R^0\ra$ with $\a_i=\dot\a_i+\frac{1}{n_i}\d_i,$ so $a=\sum_{i=1}^mr_i\dot\a_i+\sum_{i=1}^m\frac{r_i}{n_i}\d_i.$ This implies that $a=\sum_{i=1}^m\frac{r_i}{n_i}\d_i.$ Therefore we have $n_1\cdots n_m a\in\la R^0\ra.$

For the last assertion, we just need to assume $R$ is of  type $BC(1,1).$ In this case, regarding the description $R=\cup_{\dot\a\in \dot R}(\dot\a+T_{\dot\a})$ for $R$ as above,
$T_{\dot\a}\sub R^0$ for $\dot\a\in \dot R_{re}$ as in Case 1. Now suppose $R^0=\{0\},$ so  $T_{\dot\a}=\{0\}$ for $\dot\a\in \dot R_{re}.$ Suppose that $\dot R=\{0,\pm\ep_0,\pm\d_0,\pm2\ep_0,\pm2\d_0, \pm\ep_0\pm\d_0\}.$ Now if $r,s\in\{\pm1\}$ and $\d\in T_{r\ep_0+s\d_0},$ since $(r\ep_0,r\ep_0+s\d_0+\d)\neq 0,$ we get that $s\d_0+\d\in R$ and so $\d\in T_{s\d_0}=\{0\}.$ This shows that $R\sub \dot R$ and so $\v^0=\{0\}$ which in turn implies that $A^0=\{0\}.$ This completes the proof.
\qed

The following example shows that the condition   $X\neq A(\ell,\ell)$ is necessary in Lemma \ref{radical}. This is a phenomena occurring in the  super-version of root systems; more precisely, one knows that for an affine reflection system $(A,\fm,R)$ i.e., an  extended affine root supersystem with no nonsingular root, $R^0=\{0\}$ if and only if $A^0=\{0\};$
 see \cite{AYY}.

\begin{Examp}{\rm
(i) Suppose that $(\dot A, \fm, \dot R)$ is a locally finite root supersystem of type $X= A(\ell,\ell)$ for some integer $\ell\geq2$  as in Theorem \ref{classification II} with Weyl group $\w.$ Suppose that  $\sg$ is a symbol and set $A:=\dot A\op \bbbz \sg.$ Fix $\d^*\in \dot R_{ns}^\times$ and note that $-\d^*\not \in \w\d^*.$ Set $R:=\dot R_{re}\cup \pm (\w\d^*+\sg).$ Extend the form on $\dot A$ to a form on $A$ denoted again by $\fm $ such that $\sg$ is an element of the radical of  this new  form. Set  $B:=\la R\ra.$ We claim that the form $\fm$ restricted to $B$ is degenerate; indeed,  since $\dot R$ is of real type, there is $n\in\bbbz\setminus\{0\}$ such that $n\d^*\in \la \dot R_{re}\ra\sub B,$ so $n\sg=n(\d^*+\sg)-n\d^*\in B$ which in turn implies that $n\sg$ is an element of the radical of the form on $B.$ One can check that for $\a\in \rim$ and $\b\in R$ with $(\a,\b)\neq0,$ we have  either $\a+\b\in R$ or $\a-\b\in R.$
Next we note that  $R^0=\{0\},$  the root string property is satisfied for $\dot R_{re}$ and that for $\a\in \rre^\times$ and $\b\in R,$ we have $r_\a\b\in R.$   These  together  with  the same argument as in  \cite[Lem. 3.12]{you2} imply that the  root string property is satisfied for $R.$ These all together imply that  $R$ is an extended affine root supersystem with $R^0=\{0\}$ but it is not a locally finite root supersystem as the form on $B$ is degenerate.

(ii) Suppose that $(\dot A,\fm, \dot R)$ is  a locally finite root supersystem of type  $A(1,1)$ as in Theorem \ref{classification II}. Suppose that  $\sg$ is a symbol and set $A:=\dot A\op \bbbz \sg.$ Set $R:=\dot R_{re}\cup (\dot R_{ns}^\times\pm\sg).$ Extend the form on $\dot A$ to a form on $A$ denoted again by $\fm $ such that $\sg$ is  an element of the radical of  this new  form. As above, the form restricted to   $B:=\la R\ra$ is degenerate and $R$ is an extended affine root supersystem, with $R^0=\{0\},$ which is not a locally finite root supersystem.
}
\end{Examp}

\begin{lem}\label{scalar-mult}
Suppose that $(A,\fm,R)$ is a locally finite  root supersystem. Then we have the following:

(i)
There is a sub-supersystem $S$ of $R$ with $R_{ns}=S_{ns}$ and $\la R\ra=\la S\ra $ such that  for $\a\in S$ and $\d\in S_{ns}$ with $(\a,\d)\neq 0,$ there is a unique $r\in\{\pm 1\}$ such that $\a+r\d\in S.$

(ii) Identify $A$ as a subset of $\bbbf\ot_\bbbz A.$ If $\d\in \rim^\times$ and  $k\in\bbbf$ with $k\d\in R,$ then $k\in\{0,\pm1\}.$
\end{lem}
\pf
$(i)$
Without loss of generality, we assume $R$ is irreducible. If $R$ is an irreducible  locally finite root supersystem of type {\small $X\neq A(1,1), BC(T,T'),$ $ C(T,T')$ $(|T|,|T'|\geq1),$} we take $R=S.$ Next suppose  $R$ is of type  {\small$X=A(1,1), BC(T,T'),$ $  C(T,T').$} We know that $\rre=R^1\op R^2$ with
$R^1,R^2$ as following:
\medskip

\begin{tabular}{|c|c|c|}
\hline
$X$&$R^1$&$R^2$\\
\hline
$A(1,1)$&$\{0,\pm\a\}$&$\{0,\pm \b\}$\\
\hline
$BC(T,T')$& $\{\pm\ep_i,\pm\ep_i\pm\ep_j\mid i,j\in T\}$& $\{\pm\d_p,\pm\d_p\pm\d_q\mid p,q\in T'\}$\\
\hline
$C(1,T)$&$\{0,\pm\a\}$& $\{\pm\ep_i\pm\ep_j\mid i,j\in T\}$\\
\hline
$C(T,T')$&$ \{\pm\ep_i\pm\ep_j\mid i,j\in T\}$& $\{\pm\d_p\pm\d_q\mid p,q\in T'\}$\\
\hline
\end{tabular}

\bigskip
Now take $S=\rim\cup S^1\cup S^2$ where $S^1, S^2$ are considered as in the following table:
$$\begin{tabular}{|c|c|c|}
\hline
$X$&$S^1$&$S^2$\\
\hline
$A(1,1)$&$\{0,\pm\a\}$&$\{0\}$\\
\hline
$BC(T,T')$& $\{0,\pm\ep_i,\pm\ep_i\pm\ep_j\mid i,j\in T,i\neq j\}$& $R^2$\\
\hline
$C(1,T)$&$\{0\}$&  $R^2$\\
\hline
$C(T,T')$&$ \{0,\pm\ep_i\pm\ep_j\mid i,j\in T,i\neq j\}$&  $R^2$\\
\hline
\end{tabular}$$
This completes the proof.

$(ii)$
{As in the proof of  Proposition \ref{lfrss}, the form $\fm$ on $A$ induces an $\bbbf$-bilinear form on $\bbbf\ot_\bbbz  A$ which is denoted by $\fm_\bbbf$ and satisfies  $(r\ot a,s\ot b)_\bbbf=rs(a,b)$ for $r,s\in \bbbf$ and $a,b\in A.$
Now suppose $\d\in R_{ns}^\times$ and  $k\neq 0$ with $k\d\in R.$ Since $\d\in R_{ns}^\times,$ there is $\b\in R$ with $(\d,\b)\neq 0.$ Then for $\b':=r_\b(\d),$ we have $(\b',\b')=0$ and $(\b',\d)\neq0.$ So without loss of generality, we assume $(\b,\b)=0.$ Now as $(\b,\d)\neq 0$ and $(\b,k\d)=(\b,k\d)_\bbbf=k(\b,\d)\neq 0,$ there are $r,s\in\{\pm 1\}$ such that $\b+r\d,\b+sk\d\in R.$ Now we have
$$\begin{array}{ll}
(\b+sk\d,\b+sk\d)=(\b+sk\d,\b+sk\d)_\bbbf=2sk(\b,\d) \andd\\
 (\b+r\d,\b+sk\d)=(\b+r\d,\b+sk\d)_\bbbf=(sk+r)(\b,\d).
 \end{array}$$
If $k=\frac{-r}{s},$ then $k\in\{\pm 1\}$ and so we are done; otherwise, $\frac{sk+r}{sk}=\frac{2(\b+r\d,\b+sk\d)}{(\b+sk\d,\b+sk\d)}$ is an integer number. This implies that  $k\in\{\pm 1\}.$
}
\qed

\begin{Def}{\rm Suppose that $(A,\fm,R)$ is a locally finite root supersystem.  A  subset $\Pi$  of $R$ is called an {\it integral base} for $R$ if $\Pi$ is a $\bbbz$-basis for $A.$ An integral base $\Pi$ of $R$ is called a {\it base} for $R$ if  for each $\a\in \rcross,$ there are $\a_{1},\ldots,\a_{n}\in\Pi$  (not necessarily distinct) and $r_1,\ldots,r_n\in\{\pm1\}$ such that $\a=r_1\a_{1}+\cdots+r_n\a_{n}$ and  for all $1\leq t\leq n,$ $r_1\a_{1}+\cdots+r_t\a_{t}\in \rcross.$
}
\end{Def}

\begin{lem}\label{base}
Suppose that $(A,\fm,R)$ is an irreducible locally finite root supersystem of type $X.$ Then $R$ contains an integral base; in particular, $A$ is a free abelian group. Moreover, if  $X\neq A(\ell,\ell),$ $R$ possesses a base.
\end{lem}
\pf Contemplating \cite[Lem. 5.1]{LN2} and \cite[\S 10.2]{LN}, we assume that $ R_{ns}\neq \{0\}$ and take $R$ to be one of the root supersystems listed in Theorems \ref{classification I} or \ref{classification II}. In what follows for  index sets $T$ and $T'$ with $|T|,|T'|\geq 2$ and a positive integer $\ell,$ we use the following notations:

{\small $${
\begin{tabular}{|c|l||c|l|}
\hline
{\tiny$A_T$}&$\{\ep_i-\ep_j\mid i,j\in T\}$&
{\tiny$ BC_1$}&$\{0,\pm\ep_0,\pm2\ep_0\},\{0,\pm\d_0,\pm2\d_0\}$\\
\hline
{\tiny$A_{T'}$}&$\{\d_p-\d_q\mid p,q\in T'\}$& {\tiny$B_T$}&$\{0,\pm\ep_i,\pm\ep_i\pm\ep_j\mid i,j\in T,i\neq j\}$\\
\hline
{\tiny$C_T$}&{\footnotesize$\{\pm\ep_i\pm\ep_j\mid i,j\in T\}$}&{\tiny$B_{T'}$}&$\{0,\pm\d_p,\pm\d_p\pm\d_q\mid p,q\in T',p\neq q\}$\\
\hline
{\tiny$C_{T'}$}&{\footnotesize$\{\pm\d_p\pm\d_q\mid p,q\in T'\}$}&
{\tiny$A_1$}& $\{0,\pm\ep_0\},\{0,\pm \d_0\},\{0,\pm\gamma_0\}$\\
\hline
{\tiny$D_T$}&{\footnotesize$B_T\cap C_{T}$}&
{\tiny$A_\ell$}& $\{\d_i-\d_j\mid 1\leq i,j\leq \ell+1\}$\\
\hline
{\tiny$ BC_{T'}$}&{\footnotesize$B_{T'}\cup C_{T'}$}&
{\tiny$A_\ell$} &$\{\ep_i-\ep_j\mid 1\leq i, j\leq \ell+1\}$\\
\hline
{\tiny$BC_T$}&$B_T\cup C_{T}$&{\tiny$G_2$}& $\{\hbox{\footnotesize$0,\pm(\ep_i-\ep_j),\pm(2\ep_i-\ep_j-\ep_t)$}\mid\hbox{\footnotesize $\{i,j,t\}=\{1,2,3\}$}\}$\\
\hline
\end{tabular}}
$$}

 In addition, we   fix $t_0\in  T$ and $p_0\in T'$ and consider the notations as in Theorems  \ref{classification I} and \ref{classification II}. We next take  $\Pi$ to be as in the following table:
{\small $${
\begin{tabular}{|c|l|}
\hline
 \hbox{type}& \parbox{3.5in}{\begin{center}$ \Pi$\end{center}}\\
 \hline
$\dot A(0,T)$ &$\{\a^*,\ep_t-\ep_{t_0}\mid t\in T\setminus\{t_0\}\}$\\
\hline
$\dot C(0, T)$ &$\{\a^*,2\ep_{t_0},\ep_t-\ep_{t_0}\mid t\in T\setminus\{t_0\}\}$\\
\hline
$\dot A(T,T')$ &$\{\a^*,\ep_t-\ep_{t_0},\d_{t'}-\d_{p_0}\mid t\in T\setminus\{t_0\},t'\in T'\setminus\{p_0\}$\\
\hline
$A(\ell,\ell)$ &$\{\ep_i-\ep_{i+1},\omega^1_1+\omega^1_2,\d_r-\d_{r+1}\mid 1\leq i\leq \ell-1,1\leq r\leq\ell\}$\\
\hline
$B(T,T')$& $\{\ep_{t_0},\ep_t-\ep_{t_0},\d_p-\ep_{t_0}\mid t\in T\setminus\{t_0\},p\in T'\}$\\
\hline
$BC(T,T')$ &$\{\ep_{t_0},\ep_t-\ep_{t_0},\d_p-\ep_{t_0}\mid t\in T\setminus\{t_0\},p\in T'\}$\\
\hline
$BC(1,1)$ &$\{\ep_0,\ep_0+\d_0\}$\\
\hline
$BC(1,T)$ &$\{\ep_{0},\ep_{t_0},\ep_t-\ep_{t_0}\mid t\in T\setminus\{t_0\}\}$\\
\hline
$D(T,T')$ & $\{2\d_{p_0},\d_p-\d_{p_0},\ep_t-\d_{p_0}\mid p\in T'\setminus\{p_0\},t\in T\}$\\
\hline
$C(T,T')$& $\{2\ep_{t_0},\ep_t-\ep_{t_0},\d_p-\ep_{t_0}\mid t\in T\setminus\{t_0\},p\in T'\}$\\
\hline
$B(1,T)$&$\{\ep_{0},\ep_0-\ep_t\mid t\in T\}$\\
\hline
$C(1,T')$&$\{\ep_{0},\frac{1}{2}\ep_0-\d_p\mid p\in T'\}$\\
\hline
$AB(1,3)$& $\{\ep_1-\ep_2,\ep_2-\ep_3,\ep_3,\frac{1}{2}(\ep_0-\ep_1-\ep_2-\ep_3)\}$\\
\hline
$D(1,T)$& $\{\ep_0,\frac{1}{2}\ep_0-\ep_t\mid t\in T \}$\\
\hline
$B(T,1)$& $\{\ep_0,\ep_0-\ep_t\mid t\in T \}$\\
\hline
$G(1,2)$&$\{ \ep_0,\ep_0-\ep_1+\ep_2,2\ep_1-\ep_2-\ep_3\}$\\
\hline
$D(2,1,\lam)$ &$\{\ep_0,\d_0,\frac{1}{2}\ep_0+\frac{1}{2}\d_0+\frac{1}{2}\gamma_0\}$\\
\hline
$D(2,T)$ & $\{\ep_0,\d_0,\frac{1}{2}\ep_0+\frac{1}{2}\d_0+\ep_{t_0},\ep_t-\ep_{t_0}\}\mid t\in T\setminus\{t_0\}\}$\\
\hline
\end{tabular}
}$$}

{One can check that $\Pi$ is an integral base for $R$ and that if $R$ is  not of type $A(\ell,\ell),$ $\Pi$ is a base for $R.$}
\qed

\medskip
{Using the same argument as in  Lemma 3.1 of \cite{you2} and contemplating Lemma \ref{base}, one has the following lemma:}
\begin{lem}\label{full}
(i) Suppose that $\Pi$ is a base for an irreducible locally finite root supersystem $(A,\fm, R).$  Then for each  finite subset  $X\sub\Pi,$ there is a  finite subset  $Y_X\sub \Pi$ such that $X\sub Y_X$ and the form restricted to $\la Y_X\ra$ is nondegenerate. {Moreover, if $X$ is connected, we can choose $Y_X$ to be connected.}

(ii) If $\Pi$ is a connected integral base for a locally finite root supersystem $R,$ then $R$ is irreducible.

(iii) Suppose that $R$ is an infinite irreducible locally finite root supersystem in an additive abelian group $A.$  Then there is a base $\Pi$ for $R$ and a class $\{R_\gamma\mid \gamma\in\Gamma\}$ of finite irreducible $\bbbz$-linearly closed  sub-supersystems of $R$ of the same type as $R$ such that $R$ is the direct union of $R_\gamma$'s and for each $\gamma\in\Gamma,$ $\Pi\cap R_\gamma$ is a base for $R_\gamma.$  {In particular, each finite subset of $R$ lies in a finite $\bbbz$-linearly closed sub-supersystem.}

\end{lem}

\section{Structure Theorem}
In this section, we give a description of    the structure of extended affine root supersystems.
The following proposition is a generalization of Proposition 5.9 of \cite{huf} to extended affine root supersystems.

\begin{Pro}
\label{super-root}
Suppose that $A$ is an additive abelian group equipped with a  symmetric form.
Consider the induced form  $\fm^{\bar{}}$ on  $\bar A=A/A^0$ and suppose that
$\bar\;:A\longrightarrow \bar A$ is the {canonical epimorphism}. Assume that  $S$ is a subset of $A^\times:=A\setminus A^0$ and set $B:=\la S\ra.$ If
\begin{itemize}
\item $(\bar B, \fm^{\bar{}}\mid_{\bar B\times \bar B}, \bar S\cup\{0\})$ is a locally finite root supersystem,
\item {$S=-S$ and  $\a-\frac{2(\a,\b)}{(\b,\b)}\b\in S$ for $\b\in S_{re}^\times$ and  $\a\in S,$}
\item for $\a\in S_{ns}$ and $\b\in S$ with $(\a,\b)\neq 0,$  $\{\b+\a,\b-\a\}\cap S\neq\emptyset,$
\end{itemize}
 then $R:=S\cup((S-S)\cap A^0)$ is a tame extended affine root supersystem in its $\bbbz$-span.
\end{Pro}
\pf  To show that $R$ is a tame extended affine root supersystem, we just need to prove that the root string property holds. Suppose that $\a\in \rre^\times$ and $\b\in R.$


\textbf{\underline{Step 1.}} {$\b\not\in{R}_{ns}^\times:$ We know that $\bar \a,\bar \b$ are two elements of the locally finite root supersystem $\bar S\cup\{0\}.$ So using Lemma \ref{full}, there is a finite sub-supersystem $\Phi$ with $\bar \a,\bar \b\in\Phi$ such that $(\bbbz\bar \a+\bbbz\bar \b)\cap \bar R\sub \Phi.$ Set $$R_{\a,\b}:=R\cap(\bbbz\a+\bbbz\b)\andd  X_{\a,\b}:=R\cap (\b+\bbbz\a).$$ Then as  $(\bar R_{\a,\b})_{re}$ is  invariant under the reflections based on its nonzero elements, it  is a subsystem of   the finite root system $\Phi_{re}.$ Now we carry out the proof through the  following two cases:}

Case 1.
{$(\bbbz\a+\bbbz\b)\cap A^0\neq\{0\}$: In this case, $\bar\a,\bar \b$ are $\bbbz$-linearly dependent elements of $(\bar R_{\a,\b})_{re},$ so $(\bar R_{\a,\b})_{re}$ is a finite root system of rank $1,$
 in other words, it is either of type $A_1$ or $BC_1.$ Then there is $\d\in A^0$ such that
$$X_{\a,\b}\sub \left\{\begin{array}{ll}
\{0,\pm1,\pm2\}\a+\d&\hbox{if  $\bar\b\in\{0,\pm2\bar\a,\pm\bar\a\}$}\\
\{\b,\b\mp\a\}=\{\b,r_\a(\b)\} & \hbox{if $\bar\a\in\{\pm2\bar\b\}.$}
\end{array}\right.$$
So the root string property holds if $\bar\a\in\{\pm2\bar\b\}.$  If  $\bar\b\in\{0,\pm2\bar\a,\pm\bar\a\},$ then $X_{\a,\b}=Y\a+\d$ where $Y$ is a subset of $\{0,\pm1,\pm2\}.$ As $R$ is invariant under $r_\a,$ we have $Y=-Y.$ If $2\in Y,$  then $\d+\a=-r_{\d+2\a}\a\in R$ that is $1\in Y.$ Also if $1\in Y,$ than  $\d+\a\in S,$ and so  $\d=(\d+\a)-\a\in (S-S)\cap A^0\sub R,$ that is $0\in Y.$ We conclude that either $Y=\{0,\pm1,\pm 2\}$ or $Y=\{0,\pm 1\}$ and so the root string property holds.
}

Case 2. {$(\bbbz\a+\bbbz\b)\cap A^0=\{0\}$: If   $\bar \a,\bar \b$ are $\bbbz$-linearly dependent, we get the result as in Case 1. We next suppose  $\bar \a,\bar \b$ are $\bbbz$-linearly independent.  We claim  that the form restricted to $\bbbz\a+\bbbz\b$ is nondegenerate.  We suppose that $r\a+s\b$ is an element of the radical of the form on $\bbbz\a+\bbbz\b$ and prove that $r=s=0.$ If either $r=0$ or $s=0,$ we are done. So we assume  $r,s\neq 0$ and get a contradiction.  We have  $r(\a,\a)+s(\b,\a)=(r\a+s\b,\a)=0$ and $r(\a,\b)+s(\b,\b)=(r\a+s\b,\b)=0.$
This implies that $(\bar \a,\bar \b)/(\bar \a,\bar \a)=-r/s$ and $(\bar \a,\bar \b)/(\bar \b,\bar \b)=-s/r.$ But  $\bar\a,\bar \b$ are two  $\bbbz$-linearly independent roots of the finite root system $\Phi_{re},$ so we get  $4=(2r/s)(2s/r)=\frac{2(\bar \a,\bar \b)}{(\bar \b,\bar \b)}\frac{2(\bar \a,\bar \b)}{(\bar \a,\bar \a)}\in\{0,1,2,3\},$ a contradiction. Therefore,  the form restricted to $\bbbz\a+\bbbz\b$ is nondegenerate. Also it is immediate that  $R_{\a,\b}$ satisfies  (S1)-(S3) and (S5).
We next take $\phi$ to be the restriction of ``$\;\bar{\;}\;$'' to $(R_{\a,\b})_{re}.$ Since $(\bbbz\a+\bbbz\b)\cap A^0=\{0\},$ $\phi$ is an embedding into $\Phi_{re};$
in particular, $(R_{\a,\b})_{re}$ is finite. This in particular implies that the  root string property holds in $R_{\a,\b}$ (see \cite[Lem. 3.10]{you2}) and so we are done in this case.
}

\textbf{\underline{Step 2.}}
{$\b\in R_{ns}^\times:$ If $X_{\a,\b}\cap (R\setminus R_{ns}^\times)\neq\emptyset,$ then we get the result by Step 1. So we assume $X_{\a,\b}\sub  S_{ns}^\times.$  Therefore, for $k\in \bbbz,$ $k\a+\b\in R$ implies that $(k\a+\b,k\a+\b)=0.$ Since $(\a,\a)\neq 0,$ this gives $X_{\a,\b}\sub\{\b,r_\a\b\}.$ Since $r_\a\b,\b\in X_{\a,\b},$ we get $X_{\a,\b}=\{\b,r_\a\b\}.$ If $(\b,\a)=0,$ this gives $X_{\a,\b}=\{\b\},$ so the string property holds. If $(\b,\a)\neq 0,$ then $\b+\a$ or $\b-\a$ lies in $X_{\a,\b},$ so $r_\a\b$ is either $\b+\a$ or $\b-\a;$ in both cases the root string property holds. This completes the proof.}\qed

\smallskip

\medskip

In \cite[\S 3]{N1} and \cite[Thm. 1.13]{AYY},  the authors give  the  structure of an  affine reflection system i.e., an extended affine root supersystem whose set of nonsingular roots is $\{0\}.$
In  the following theorem, we give the structure of extended affine root supersystems. We  see that the notion of extended affine root supersystems is in fact a generalized notion of root systems extended by an abelian group introduced by Y. Yoshii \cite{yos2}. More precisely, we show that associate to each extended affine root supersystem  $(A,\fm,R)$ of type $X,$ there is  a locally finite root supersystem $\dot R$ as well as  a class $\{S_{\dot\a}\}_{\dot\a\in \dot R}$ of subsets of $A^0$ such that $R=\cup_{\dot\a\in\dot R}(\dot\a+ S_{\dot\a}).$ If $X\neq  A(\ell,\ell),C(1,2),C(T,2),BC(1,1),$ then the interaction of $S_{\dot\a}$'s
 results in a nice characterization of $R.$

In what follows by a {\it reflectable set}   for a locally finite root system $S,$ we mean a subset $\Pi$ of $S\setminus\{0\}$ such that $W_\Pi(\Pi)$ coincides with the set of nonzero reduced roots $S^\times_{red}=S\setminus\{2\a\mid\a\in S\},$ in which  $W_\Pi,$ is  the subgroup of the Weyl group generated by $r_\a$ for all $\a\in \Pi;$ see \cite{AYY}.
 We also  recall from \cite{L} that a {\it symmetric reflection subspace} (or s.r.s for short) of an additive abelian group $A$ is
a nonempty subset $X$ of $A$ satisfying  $X -2X\sub X;$ we mention that  a symmetric reflection subspace satisfies $X=-X.$ A symmetric reflection subspace $X$ of an additive abelian group $A$ is called a {\it pointed reflection subspace} (or p.r.s for short) if $0\in X.$ Before stating  the structure theorem of extended affine root supersystems, we make a convention that if $\dot R$ is a locally finite root supersystem with  decomposition $\dot R_{re}=\op_{i=1}^n\dot R_{re}^i$ of $\dot R_{re}$ into irreducible subsystems, by $\dot R_{*},$ $*=sh,lg,ex,$ we mean $\cup_{i=1}^n (\dot R^i_{re})_{*}.$

\begin{The} Suppose that $(\dot A,(\cdot,\cdot\dot),\dot R)$ is an irreducible locally finite root supersystem of type $X$ with $\dot R_{ns}\neq\{0\},$ as in Theorems \ref{classification I} and \ref{classification II}, and $A^0$ is an additive  abelian group. Extend the form $(\cdot,\cdot\dot)$ to the  form $\fm$ on $\dot A\op A^0$ whose radical is $A^0.$

(i) Suppose that {\small $X\neq A(\ell,\ell),  BC(T,T'),C(T,T'),C(1,T)$,} $F$ is a subgroup of $A^0$ and $S$ is a pointed reflection subspace of $A^0$ such that
$$\begin{array}{l}\la S\ra= A^0, \;\; F+S\sub S,\;\; 2S+F\sub F \andd \\
S=F\; \hbox{if }X\neq B(T,T'),B(T,1),B(1,T).\end{array}$$  Then $$R:=(S-S)\cup (\dot R_{sh}+S)\cup ((\dot R^\times\setminus \dot R_{sh})+F)$$ is a tame irreducible extended affine root supersystem of type $X.$ Conversely, each tame irreducible extended affine root supersystem of type $X$  arises in this manner.

(ii) Suppose that  {\small $X=BC(1,T),BC(T,T'),$ $|T|,|T'|>1,$} $F$ is a subgroup of $A^0,$  $S$ is a pointed reflection subspace of $A^0$ and $E_1,E_2$ are two symmetric reflection subspaces of $A^0$ such that
$$\begin{array}{l}\la S\ra= A^0,\; \parbox{4in}{$\{\sg+\tau,\sg-\tau\}\cap (E_1\cup E_2)\neq\emptyset;\;\;$ $\sg,\tau\in F,$}\\\\
F+S\sub S,\;\; 2S+F\sub F,\; 2F+E_i\sub E_i\; (\hbox{if } (\dot R_{re}^i)_{lg}\neq\emptyset),\;  F+E_i\sub F,\\ S+E_i\sub S,\; E_i+4S\sub E_i\; (i=1,2).
\end{array}$$  Then $$R:=(S-S)\cup (\dot R_{sh}+S)\cup (\dot R^1_{ex}+E_1)\cup (\dot R^2_{ex}+E_2)\cup ((\dot R_{lg}\cup\dot R^\times_{ns})+F)$$ is a tame irreducible extended affine root supersystem of type $X;$ conversely each tame irreducible extended affine root supersystem of type $X$ arises in this manner.

(iii) Suppose that  {\small $X=C(1,T'),$ $|T'|>2,$} $F$ is a subgroup of $A^0,$ $S$ is a pointed reflection subspace of $A^0$ and $L$ a symmetric reflection subspace of $A^0$ such that
$$\begin{array}{l}F= A^0,\;F=S\cup L,  L+2F\sub L.
\end{array}$$  Then $$R=F\cup (\dot R^1_{sh}+S)\cup ((\dot R^2_{sh}\cup\dot R^\times_{ns})+F)\cup (\dot R^2_{lg}+L)$$ is a tame irreducible extended affine root supersystem of type $X.$ Conversely, each tame irreducible extended affine root supersystem of type $C(1,T'),$ $|T'|>2,$ arises in this manner.

(iv)  Suppose that  {\small $X=C(T,T'),$ $|T|\geq2,|T'|>2,$} $F$ is a subgroup of $A^0,$ $L_1$ is a pointed reflection subspace of $A^0$ and $L_2$ is a symmetric reflection subspace of $A^0$ such that
$$\begin{array}{l}F= A^0,\;F=L_1\cup L_2,\; L_i+2F\sub L_i\;\;\; (i=1,2).
\end{array}$$  Then\footnote{I thank Gholamreza Behboodi for  pointing   out that one can write the conditions in ($iii$) and ($iv$) in this simple form. In the published version of the paper, these conditions have been written respectively as
$$\{\sg+\tau,\sg-\tau\}\cap (S\cup L)\neq\emptyset\;(\sg,\tau\in F),
F+S\sub F,\;\; L+F\sub F,\;\; L+2F\sub L$$
and
$$\{\sg+\tau,\sg-\tau\}\cap (L_1\cup L_2)\neq\emptyset\;(\sg,\tau\in F),
L_i+F\sub F, \; L_i+2F\sub L_i\;\;\; (i=1,2).$$ } $$R=F\cup ((\dot R_{sh}\cup\dot R^\times_{ns})+F)\cup ((\dot R^1_{re})_{lg}+L_1)\cup ((\dot R_{re}^2)_{lg}+L_2)$$ is a tame irreducible extended affine root supersystem of type $X.$ Conversely, each tame irreducible extended affine root supersystem of type $C(T,T'),$ $|T|\geq2,|T'|>2,$ arises in this manner.
\end{The}
\pf  Suppose that  $(A,\fm,R)$ is  a tame irreducible extended affine root supersystem of type $X$ with $\rim\neq \{0\},$ then  by Proposition
\ref{lfrss}($iii$), $(\bar A,\fm^{\bar{}},\bar R)$ is a locally finite root supersystem.
Fix a  subset $\dot \Pi$ of $ R$ such that $\bar{\dot\Pi}$ is the corresponding  $\bbbz$-basis for $\bar A$ as introduced in Lemma  \ref{base}.
Take $\dot A:=\la\dot\Pi\ra$ as well as $\dot R:=\{\dot a\in \dot A\mid \exists \eta\in A^0;\;\dot a+\eta\in R\}.$ One can see that $A=\dot A\op A^0,$
%
and  that $\dot R$ is a locally finite root supersystem in $\dot A$ isomorphic to $\bar R.$
Without loss of generality, by multiplying the form  $\fm$ to a nonzero scalar, we may assume $\dot R$ is one of the locally finite root supersystems as in Theorems \ref{classification I} and \ref{classification II}  with the decomposition  $\dot R_{re}=\op_{i=1}^n\dot R_{re}^i$ for $\dot R_{re}$ into irreducible subsystems and that $\dot \Pi$ is as in Lemma \ref{base}.

\textbf{Claim 1.} If $\dot R$ is of type $X\neq A(\ell,\ell), C(T,T'),C(1,T),$ then $R$ contains a reflectable set for $\dot R_{re}:$  We note that if $\dot R$ is of imaginary type, then $\dot \Pi\cap \dot R_{re}\subset R$ is a reflectable set for $\dot R_{re}.$ So we suppose that $\dot R$ is of real type and  carry out the proof through the following cases:
\begin{itemize}
\item{Case 1.} \underline{$\dot R$ is  of type  $AB(1,3)$}:
In this case, since $\ep_1-\ep_2,\ep_2-\ep_3,\ep_3\in\dot R\cap R,$ we get that $\ep_1+\ep_2=r_{\ep_2-\ep_3}r_{\ep_2-\ep_1}r_{\ep_3}(\ep_2-\ep_3)\in\dot R\cap R.$ We note that for  $\dot \a:=\ep_1+\ep_2,\dot \b:=\ep_3,\dot\gamma:=\frac{1}{2}(\ep_0-\ep_1-\ep_2-\ep_3)\in R\cap \dot R,$ we have $\dot\a-\dot\gamma,\dot\b-\dot\gamma\not\in \dot R.$ Therefore
$\bar{\dot\a}-\bar{\dot\gamma},\bar{\dot\b}-\bar{\dot\gamma}\not \in \bar R$
and so $\dot\a-\dot\gamma,\dot\b-\dot\gamma\not\in  R.$ This together with the fact that $(\dot\a,\dot\gamma)\neq0$  and $(\dot\b,\dot\gamma)\neq0,$ implies that  $\dot\a+\dot\gamma,\dot\b+\dot\gamma\in R,$ and so $\dot\eta:=\dot\a+\dot\gamma,\dot\zeta:=\dot\b+\dot\gamma\in  \dot A\cap R\sub R\cap \dot R.$
Again as $(\dot\eta,\dot\zeta)\neq0,$ the same argument as above implies that $\ep_0=\dot\eta+\dot\zeta\in R\cap \dot R.$ So we are done as  $\{\ep_0,\ep_3,\ep_1-\ep_2,\ep_2-\ep_3\}$ is a reflectable set for $\dot R_{re}.$

\item{Case 2.} \underline{$\dot R$ is  of type  $D(2,1,\lam)$}:
We know that   $\dot\eta:=\frac{1}{2}\ep_0+\frac{1}{2}\d_0+\frac{1}{2}\gamma_0,\ep_0,\d_0\in R\cap \dot R.$ Since $ \dot\eta+ \ep_0, \dot\eta+\d_0\not\in  \dot R,$ we get that  $\bar{\dot \eta}+\bar \ep_0,\bar {\dot\eta}+\bar \d_0\not\in \bar R$ and so $\dot \eta+ \ep_0, \dot\eta+ \d_0\not\in  R.$ Also as $(\dot\eta,\ep_0)\neq0$ and $(\dot\eta,\d_0)\neq0,$ we have $\dot\zeta:=\dot\eta-\ep_0=-\frac{1}{2}\ep_0+\frac{1}{2}\d_0+\frac{1}{2}\gamma_0,\dot\xi:=\dot\eta-\d_0=\frac{1}{2}\ep_0-\frac{1}{2}\d_0+\frac{1}{2}\gamma_0\in R\cap \dot A.$ Again as $(\dot\zeta,\dot\xi)\neq0$ and $\dot \xi-\dot \zeta\not\in \dot R,$ we get that $\gamma_0= \dot\xi+\dot\zeta\in \dot A\cap R $ and so $\gamma_0\in \dot R\cap R.$ Therefore, $\{\ep_0,\d_0,\gamma_0\},$ which is a reflectable  set  for $\dot R_{re},$ is contained in $R.$

\item{Case 3.} \underline{$\dot R$ is  of type  $D(2,T)$}:
Using the same argument as in the previous case, we get that $2\ep_{t_0}\in \dot R\cap R$ and so  $\{\ep_0,\d_0,2\ep_{t_0},\ep_t-\ep_{t_0}\mid t\in T\setminus\{t_0\}\},$ which is a reflectable set for $\dot R_{re},$ is contained in $R.$

\item{Case 4.} \underline{$\dot R$ is  of type  $D(T,T')$}: Since for $i\in T,$ $\ep_i-\d_{p_0},2\d_{p_0}\in R\cap\dot R$ and
$\ep_i-\d_{p_0}-2\d_{p_0}\not\in \dot R,$ as above one  concludes that $\ep_i+\d_{p_0}\in R\cap \dot R.$ Now for $i,j\in T$ with $i\neq j,$ we have $\dot \a:=\ep_i-\d_{p_0},\dot\b:=\ep_i+\d_{p_0},\dot\gamma:=\ep_j-\d_{p_0}\in R\cap \dot R$ with $(\dot\a,\dot\gamma)\neq 0$ and $(\dot\b,\dot\gamma)\neq0,$ but $\dot\a+\dot\gamma,\dot\b-\dot\gamma\not\in \dot R,$ so as above, we get that $\ep_i+\ep_j,\ep_i-\ep_j\in R \cap \dot R.$ This completes the proof in this case as  $\{2\d_{p_0},\d_p-\d_{p_0},\ep_i\pm\ep_j\mid p\in T'\setminus\{p_0\},i\neq j\in T\}\sub R$ is a reflectable set for $\dot R_{re}.$

\item{Case 5.} \underline{$\dot R$ is  of type  $ D(1,T)$}: We know that for $t\in T,$ $\ep_0,\frac{1}{2}\ep_0-\ep_t\in\dot R\cap R$ and that $\ep_0+(\frac{1}{2}\ep_0-\ep_t)\not\in \dot R,$ so as before, we get that $\frac{1}{2}\ep_0+\ep_t\in R\cap\dot R.$ Using the same argument as in  the previous case, we get that $\ep_r\pm\ep_s\in R\cap \dot R$ for all $r,s\in T$ with $r\neq s.$ This completes the proof in this case.

\item{Case 6.} \underline{$\dot R$ is of type  $B(T,T'), BC(T,T'), B(1,T), B(T,1), G(1,2)$}:  In these cases, for $\dot\Pi_{re}:=\dot\Pi\cap
\dot R_{re}$ and $ \dot\Pi_{ns}:=\dot\Pi\cap \dot R_{ns},$ the set $\dot\Pi_{re}\cup ((\dot\Pi_{ns}\pm\dot\Pi)\cap (\dot R_{re})_{red}^\times),$ which  is  (as above) a subset of $R,$ is a reflectable set for $\dot R_{re}.$
\end{itemize}

\textbf{Claim 2.} If $X=C(T,T'),C(1,T'),$ then $\dot R_{re}^\times\setminus (\dot R_{re}^2)_{lg}\sub R:$  We know that $\dot \Pi\sub R \cap \dot R.$ So as in Case 5 of the the proof of Claim 1, we get that $\pm\d_p\pm\d_q\in R_{re}$ for all $p,q\in T'$ with $p\neq q.$ Moreover, if $X=C(T,T'),$ then  for $t\in T\setminus\{t_0\},$ since $2\ep_{t_0},\ep_t-\ep_{t_0}\in R\cap \dot R,$   we have that $$\begin{array}{ll}\ep_t+\ep_{t_0}=r_{2\ep_{t_0}}(\ep_t-\ep_{t_0})\in\dot R\cap R,&\ep_r-\ep_t=r_{\ep_t-\ep_{t_0}}(\ep_r-\ep_{t_0})\in\dot R\cap R,\\ \ep_r+\ep_t=r_{\ep_t-\ep_{t_0}}(\ep_r+\ep_{t_0})\in\dot R\cap R,&2\ep_t=r_{\ep_{t_0}-\ep_t}(2\ep_{t_0})\in\dot R\cap R,\\
\end{array}$$
 for all $r,t\in T$ with $r, t\in T\setminus\{t_0\}$ and $r\neq t.$ This completes the proof of the claim.

\textbf{Claim 3.} For $\dot\a\in \dot R,$ set $$S_{\dot\a}:=\{\eta\in A^0\mid \dot\a+\eta\in R\}.$$ Then we have
\begin{equation}\label{final}\left\{\begin{array}{ll}
0\in S_{\dot\a}& \hbox{if {\small $X\neq A(\ell,\ell),C(T,T'),C(1,T')$ $\;\&\;$ $\dot\a\in (\dot R_{re})_{red},$}}\\
0\in S_{\dot\a}& \hbox{if {\small $X=C(T,T'),C(1,T')$ $\;\&\;$ $\dot\a\in\dot R_{re}^\times\setminus(\dot R_{re}^2)_{lg},$}}
\end{array}\right.\end{equation} and that

\begin{equation}\label{**}
\parbox{3in}{  if $X\neq A(\ell,\ell), C(1,2),C(T,2),$  $S_{\dot\a}=S_{\dot\b}$ for all $\dot\a,\dot\b\in \dot R^i_{re}\setminus\{0\}$\;($1\leq i\leq n$) with $(\dot\a,\dot\a)=(\dot\b,\dot\b):$}
\end{equation}
To show this, suppose that $\dot\a\in\dot R_{re}^\times,\dot\b\in \dot R,$  $\eta\in S_{\dot\a}$ and $\zeta\in S_{\dot\b},$ then \begin{eqnarray*}r_{\dot\a+\eta}(\dot\b+\zeta)=\dot\b+\zeta-\frac{2(\dot\a,\dot\b)}{(\dot\a,\dot\a)}(\dot\a+\eta)=
r_{\dot\a}(\dot\b)+(\zeta-\frac{2(\dot\a,\dot\b)}{(\dot\a,\dot\a)}\eta).
\end{eqnarray*}
This means that \begin{equation}
\label{rel-s-alpha}
S_{\dot\b}-\frac{2(\dot\a,\dot\b)}{(\dot\a,\dot\a)}S_{\dot\a}\sub S_{r_{\dot\a}(\dot\b)};\;\;\;\;\;(\dot\a\in\dot R_{re}^\times,\;\dot\b\in \dot R).
\end{equation}
This in turn implies that  \begin{equation}\label{****}
S_{\dot\b}-2S_{\dot\b}\sub S_{-\dot\b};\;\;\;\;(\dot\b\in\dot R_{re}^\times).
\end{equation}

Now suppose that $\dot R$ is of type $X=C(T,T'),C(1,T').$ Using Claim 2, we have $0\in S_{\dot\a}$ for $\dot\a\in \dot R_{re}^\times\setminus(\dot R_{re}^2)_{lg}.$
We next mention that  for a locally finite root system of type $C$ with rank greater that 2, roots of the same length are conjugate under the subgroup of the Weyl group generated by the reflection based on the short roots, therefore in this case, by (\ref{rel-s-alpha}), we get $S_{\dot\a}=S_{\dot\b}$ for all $\dot\a,\dot\b\in \dot R^i_{re}\setminus\{0\}$\;($1\leq i\leq n$) with $(\dot\a,\dot\a)=(\dot\b,\dot\b).$

Next suppose that  $X\neq A(\ell,\ell), C(T,T'),C(1,T').$ Using   Claim 1, we get that $R$ contains  a reflectable set for $\dot R_{re},$ say $\dot B.$ Now for $\dot\a\in (\dot R_{re})_{red}\setminus\{0\},$ there are $\dot\a_1,\ldots,\dot\a_{t+1}\in \dot B\sub R\cap \dot R$ such that $r_{\dot\a_1}\cdots r_{\dot\a_t}(\dot\a_{t+1})=\dot\a,$ so as $R$ and $\dot R$ are closed under the reflection actions, we get that $\dot\a\in R\cap \dot R;$ in particular $0\in S_{\dot\a}$ for $\dot\a\in (\dot R_{re})_{red}.$  These all together with (\ref{rel-s-alpha}) and the fact that for a locally finite root system, the roots of the same length are conjugate under the Weyl group action complete the proof.

{\bf Claim 4.}   Suppose that  $X\neq A(\ell,\ell), C(1,2),C(T,2),BC(1,1).$ Fix  a nonzero $\dot\d^*\in \dot R_{ns}\cap\dot \Pi\sub\dot R\cap R.$ Consider (\ref{**}) and set $$F:=S_{\dot\d^*}\andd \left\{\begin{array}{ll}
  S_i:=S_{\dot\a} & \dot\a\in (\dot R^i_{re})_{sh}\\
    L_i:=S_{\dot\a} & \dot\a\in (\dot R^i_{re})_{lg}\\
      E_i:=S_{\dot\a} & \dot\a\in (\dot R^i_{re})_{ex}
\end{array}\right.$$ for $1\leq i\leq n.$
Then
\begin{equation}\label{8}\begin{array}{ll}
\hbox{$S_i$ is a p.r.s. of $A^0$ }&\\
\hbox{$E_i$ is a s.r.s. of $A^0$}  & \hbox{if $(\dot R^i_{re})_{ex}\neq\emptyset,$  }\\
\hbox{$L_i$  is a p.r.s. of $A^0$ }& \hbox{if {\small $X\neq C(1,T'),C(T,T')$} and  $(\dot R^i_{re})_{lg}\neq\emptyset,$}\\
\hbox{$ L_2$ is a s.r.s. of $A^0$}&\hbox{if   {\small $X= C(1,T'),C(T,T');$ $|T'|>2,$} } \\
\hbox{$L_1$ is a p.r.s. of $A^0$} &\hbox{if {\small $X=C(T,T');$ $|T'|>2$}}
\end{array}\end{equation} and \begin{equation}\label{triple}
\begin{array}{ll}
(a)&S_i + L_i\subseteq S_i,\; L_i +\rho_i S_i \subseteq L_i\;\;\;\;\hbox{if }(\dot R^i_{re})_{lg}\not=\emptyset,\vspace{2mm}\\
(b)&S_i + E_i \subseteq S_i,\; E_i + 4S_i \subseteq E_i\;\;\;\;\hbox{if }\dot R^i_{re}=BC_1,\vspace{2mm}\\
(c)&L_i + E_i \subseteq L_i,\; E_i + 2L_i\subseteq E_i \;\;\;\;\hbox{if }\dot R_{re}^i=BC_P\;\;(|P|\geq 2),
\end{array}\end{equation}
in which $$
\rho_i:=(\dot\b,\dot\b)/(\dot\a,\dot\a),\qquad (\dot\a\in(\dot R_{re}^i)_{sh},\;\dot\b\in(\dot R_{re}^i)_{lg}\;\;\hbox{if $(\dot R_{re}^i)_{lg})\neq\emptyset:$}
$$
We immediately get (\ref{8}) using    (\ref{final}), (\ref{**}) and (\ref{****}).
Now   if $(\dot R_{re}^i)_{lg}\neq\emptyset,$ there are $\dot\a\in(\dot R_{re}^i)_{lg} $ and $\dot\b\in (\dot R_{re}^i)_{sh}$ such that $2(\dot\a,\dot\b)/(\dot\b,\dot\b)=\rho_i$ and $2(\dot\b,\dot\a)/(\dot\a,\dot\a)=1,$  so using (\ref{rel-s-alpha})  and (\ref{**}), we get (\ref{triple})($a$). If $\dot R_{re}^i$ is of type $BC_P$ for some index set $P$ with $|P|\geq2,$   one finds $\dot\a\in(\dot R_{re}^i)_{ex} $ and $\dot\b\in (\dot R_{re}^i)_{lg}$ such that $2(\dot\a,\dot\b)/(\dot\b,\dot\b)=2$ and $2(\dot\b,\dot\a)/(\dot\a,\dot\a)=1$  and so we get (\ref{triple})($c$). We similarly have (\ref{triple})($b$) as well.

{\bf Claim 5.} $F$ is a subgroup of $A^0,$ \begin{equation}\label{88}
F+2S_i\sub F;\;\; (1\leq i\leq n).
\end{equation}
 and for each $\dot \d\in \dot R_{ns}^\times,$ we have $S_{\dot\d}= F.$
Also
\begin{equation}\label{888}
F= \left\{\begin{array}{ll} S_i & \hbox{\small if $X\neq C(1,T'),BC(T,T'),B(T,T'),B(T,1),B(1,T);\;\;$ $ 1\leq i\leq n,$}\\
L_i& \hbox{\small if $X\neq C(T,T'),C(1,T)$ $\;\&\;$ $(\dot R^i_{re})_{lg}\neq\emptyset;\;\;\; 1\leq i\leq n,$}\\
 E_i & \hbox{\small if $X\not=BC(T,T')$ $\;\&\;$  $(\dot R^i_{re})_{ex}\neq\emptyset;\;\;\; 1\leq i\leq n,$ }\\
  S_2 & \hbox{\small if $X= C(1,T'):$}\\
\end{array}\right.\end{equation}
We know that  $\dot\w,$ the Weyl group of $\dot R,$ is generated by the reflections based on nonzero  reduced roots and that if $X=C(T,T'),C(1,T'),$ $|T'|>2,$ nonsingular roots are conjugate with $\dot\d^*$ under the subgroup of $\dot\w$ generated by the reflections based on the elements of  $(\dot R_{re}^1)_{sh}\cup(\dot R_{re}^2)_{sh}.$ So (\ref{rel-s-alpha}) together with  (\ref{final}) and the fact that  $\a\in R$ if and only if $-\a\in R,$ implies that \begin{equation}
 \left\{\begin{array}{ll}
S_{\pm\dot w\dot\d^*}= S_{\pm\dot\d^*}=\pm F&\dot w\in\dot\w \\
 0\in S_{\dot\a} &\dot \a\in \dot R_{ns}.
 \end{array} \right.
\end{equation}
Also one can easily see that
\begin{equation}\label{com4}\hbox{\footnotesize$(\dot\w\dot\d^*-\dot\w\dot\d^*)\cap \dot R^\times$}=\left\{\begin{array}{ll}
\hbox{\footnotesize$\dot R_{re}\setminus ((\dot R_1)_{sh}\cup (\dot R_2)_{sh})$}&\hbox{\tiny $X=B(T,T'),BC(T,T'), B(1,T),B(T,1),$}\\
\hbox{\small$(\dot R^1_{re})_{ex}\cup\dot R^2_{re}\setminus\{0\}$}& \hbox{\tiny $X=G(1,2),$}\\
\hbox{\small$\dot R_{re}^\times$}& \hbox{\tiny otherwise.}\\
\end{array}
\right.
\end{equation}
Moreover, if $1\leq i\leq n,$ we have
\begin{equation}\label{com5}(\dot \w\dot\d^*+ (\dot R^i_{re})_{sh})\cap \dot R^\times=\left\{
\begin{array}{ll}
(\dot R^j_{re})_{sh}&\hbox{\tiny $X=B(T,T'),BC(T,T'),B(1,T),B(T,1),$}\\
(\dot R_{re}^2)_{sh}&\hbox{\tiny $X=G(1,2),$ $i=1,$}\\
\dot \w\dot\d^*\cup (\dot R_{re}^1)_{sh}&\hbox{\tiny $X=G(1,2),$ $i=2,$}\\
\dot \w\dot\d^*& \hbox{\tiny otherwise}\\
\end{array}\right.
\end{equation}
and
\begin{equation}\label{com6}
\begin{array}{l}(\dot \w\dot\d^* + (\dot R^i_{re})_{lg})\cap \dot R^\times=\dot \w\dot\d^*;\;\; \hbox{if $(\dot R^i_{re})_{lg}\neq \emptyset,$ }\\
(\dot \w\dot\d^* + (\dot R^i_{re})_{ex})\cap \dot R^\times=\dot \w\dot\d^*;\;\; \hbox{if $(\dot R^i_{re})_{ex}\neq \emptyset.$ }
\end{array}
\end{equation}
We next note that
\begin{center}\parbox{4.1in}{if {\small $X\neq A(\ell,\ell),BC(T,T'),C(T,T'),C(1,T)$} and $\dot\a\in\dot R_{ns},\dot\b\in\dot R$ with $(\dot\a,\dot\b)\neq0,$ then there is a unique $r_{\dot\a,\dot\b}\in \{\pm1\}$ with $\dot\a+r_{\dot\a,\dot\b}\dot\b\in \dot R$ }
\end{center}
and that\begin{center}\parbox{4in}{if $X= BC(T,T'),C(T,T'),C(1,T)$ and  $\dot\a\in\dot R_{ns},\dot\b\in\dot R_{re}$ with $(\dot\a,\dot\b)\neq0,$ then there is a unique $s_{\dot\a,\dot\b}\in \{\pm1\}$ with $\dot\a+s_{\dot\a,\dot\b}\dot\b\in \dot R.$ Moreover,
\begin{equation}\label{***}
\begin{array}{l}
\hbox{ for  $\dot\b,\dot\gamma\in \dot R_{ns}$ with  $\dot\b+\dot\gamma,\dot\b-\dot\gamma\in \dot R,$ we have}\\
\hbox{$\dot\b+\dot\gamma,\dot\b-\dot\gamma\in\left\{\begin{array}{ll} \dot R_{ex}&  \hbox{if $X=BC(T,T'),$}\\ \dot R_{lg}&  \hbox{if $X=C(T,T'),$}\\
(\dot R^1_{re})_{sh}\cup (\dot R^2_{re})_{lg}&  \hbox{if $X=C(1,T)$};\end{array}\right.$}
\end{array}
\end{equation}}
\end{center}

(see Lemmas \ref{nonzeero-im} and \ref{scalar-mult}). Therefore, $$\begin{array}{ll}
S_{\dot\a}+r_{\dot\a,\dot\b}S_{\dot\b}\sub S_{\dot\a+r_{\dot\a,\dot\b}\dot\b}&(\dot\a\in\dot R_{ns},\dot\b\in\dot R,(\dot\a,\dot\b)\neq0)\\
&\hbox{\small $X\neq BC(T,T'),C(T,T'),C(1,T)$}
\end{array}$$
and
\begin{equation}\begin{array}{ll}\label{com7}
S_{\dot\a}+s_{\dot\a,\dot\b}S_{\dot\b}\sub S_{\dot\a+s_{\dot\a,\dot\b}\dot\b}&(\dot\a\in\dot R_{ns},\dot\b\in\dot R_{re},(\dot\a,\dot\b)\neq0)\\
&\hbox{\small $X= BC(T,T'),C(T,T'),C(1,T).$}
\end{array}\end{equation}
Now we drew the attention of the readers to the point that if $\dot\a,\dot\b\in \dot R_{ns}^\times$ with $\dot\a+\dot\b,\dot\a-\dot\b\in\dot R,$ although for $\sg\in S_{\dot\a},\tau\in S_{\dot \b},$  there is $r\in \{\pm1\}$ with  $(\dot\a+\sg)+r(\dot\b+\tau)\in R,$ we cannot conclude that
both  $(\dot\a+\sg)+(\dot\b+\tau)$ and $(\dot\a+\sg)-(\dot\b+\tau)$ are elements of $R.$
Considering this together with (\ref{***}) and  using (\ref{com4}), we have for $X\neq A(\ell,\ell),C(1,2),C(T,2)$ that

\begin{equation}\label{com'}F-F\sub \left\{\begin{array}{ll} E_i  & \hbox{\small if $X= B(1,T),B(T,T'),B(T,1),$  $(\dot R_{re}^i)_{ex}\neq\emptyset$    ,}\\
L_i  & \hbox{\small if $X= B(1,T),B(T,T'),B(T,1),$  $(\dot R_{re}^i)_{lg}\neq\emptyset$    ,}\\
L_i & \hbox{\small if $X= BC(T,T')$ and    $(\dot R_{re}^i)_{lg}\neq\emptyset$    ,}\\
 S_i & \hbox{\small if $X=C(T,T'),$ $|T'|>2$ ,}\\
  S_2 & \hbox{\small if $X= C(1,T'),$ $|T'|>2,$}\\
E_1,S_2,L_2  & \hbox{if $X=G(1,2),$}\\
S_i  & \hbox{\small if $X=$Remain  types under consideration, }\\
L_i  & \hbox{\small if $X=$Remain types under consideration, $(\dot R_{re}^i)_{lg}\neq\emptyset$   ,}\\
E_i & \hbox{\small if $X=$Remain types under consideration, $(\dot R_{re}^i)_{ex}\neq\emptyset$}
\end{array}\right.\end{equation}
for $1\leq i\leq n.$ Also by (\ref{com5}), we have
\begin{equation}\label{com9}F+S_i\sub\left\{
\begin{array}{ll}
S_j&\hbox{\tiny if $X=B(T,T'),BC(T,T'),B(1,T), B(T,1), G(1,2),\;\{i,j\}=\{1,2\},$}\\
F& \hbox{\tiny if $X\neq B(T,T'),BC(T,T'),B(1,T), B(T,1),G(1,2),$\;$1\leq i\leq n,$}\\
F& \hbox{\tiny if $X=G(1,2),$\; $i=2.$}
\end{array}\right.\end{equation}
In addition,  by (\ref{com6}) and (\ref{com7}), we have  \begin{equation}\label{*}F+L_i\sub F \;(\hbox{if $(\dot R^i_{re})_{lg}\neq\emptyset$)} \andd F+E_i\sub F \;(\hbox{if $(\dot R^i_{re})_{ex}\neq\emptyset$)}.\end{equation}In particular, since $0\in F,$  (\ref{com'}) imply that  \begin{equation}\label{com3}
F=\left\{\begin{array}{ll}L_i&\hbox{if $(\dot R^i_{re})_{lg}\neq\emptyset$ and  $X\neq C(T,T'),C(1,T)$}\\
E_i&\hbox{if $(\dot R^i_{re})_{ex}\neq\emptyset$, $X\neq BC(T,T')$}.\end{array}\right.
\end{equation}
Moreover, (\ref{com9}) implies that  \begin{equation}\label{equal}\parbox{4.2in}{if {\small $X=B(T,T'),BC(T,T'),B(1,T), B(T,1), G(1,2),$} then  $S_1=S_2,$ }\end{equation} so for $1\leq i\leq n,$ using  (\ref{com9}), we have \begin{equation}\label{com}F+S_i\sub\left\{
\begin{array}{ll}
S_i&\hbox{\tiny if $X=B(T,T'),BC(T,T'),B(1,T), B(T,1), G(1,2),$}\\
F& \hbox{\tiny if $X\neq B(T,T'),BC(T,T'),B(1,T), B(T,1).$}
\end{array}\right.\end{equation}
We also have using (\ref{8}) and (\ref{triple}) that  $L_i\sub S_i$ if    $(\dot R^i_{re})_{lg}\neq\emptyset,$ $E_i\sub L_i$  if  $(\dot R^i_{re})_{ex}\neq\emptyset$  and   $E_i\sub S_i$ if $\dot R_{re}^i$ is of type $BC_1.$  So by (\ref{com'}) and (\ref{equal}), we have
%
%
%

\begin{equation}\label{subtract}F-F\sub \left\{\begin{array}{ll} S_i & \hbox{\small if $X\neq C(1,T),BC(1,1);$  $ 1\leq i\leq n,$}\\
L_i& \hbox{\small if $X\neq C(T,T'),C(1,T)$ and $(\dot R^i_{re})_{lg}\neq\emptyset; 1\leq i\leq n,$}\\
 E_i & \hbox{\small if $X\not=BC(T,T')$ and  $(\dot R^i_{re})_{ex}\neq\emptyset; 1\leq i\leq n,$ }\\
  S_2 & \hbox{\small if $X= C(1,T'),$ $|T'|>2.$ }\\
\end{array}\right.\end{equation}
Therefore, we have
$$\left\{\begin{array}{ll}
F-F\stackrel{(\ref{subtract})}{\sub} S_i\stackrel{(\ref{com})}{\sub} F &  \hbox{\footnotesize $X\neq BC(T,T'), B(T,T'),B(1,T), $}\\
&\hbox{\footnotesize $ B(T,1),C(1,T'),$\; $1\leq i\leq n,$}\\\\
F-F\stackrel{(\ref{subtract})}{\sub} L_2\stackrel{(\ref{com3})}{\sub} F &\hbox{\footnotesize $X=BC(1,T), BC(T,T'), B(T,T'),B(1,T)$}\\
&\hbox{\footnotesize $ B(T,1);\;|T|,|T'|\geq 2,$}\\\\
F-F\stackrel{(\ref{subtract})}{\sub} S_2\stackrel{(\ref{com})}{\sub} F& \hbox{\footnotesize $X=C(1,T');\; |T'|>2.$}
\end{array}\right.$$
This means that for types $X\neq A(\ell,\ell), C(1,2),C(T,2),BC(1,1),$ $F-F\sub F$ and so $F$ is a subgroup of $A^0$ as   $0\in F.$
Also we get using (\ref{com}), (\ref{com3}) and (\ref{subtract}) that $$F= \left\{\begin{array}{ll} S_i & \hbox{\small if $X\neq C(1,T),BC(T,T'),B(T,T'),B(1,T),B(T,1);$  $ 1\leq i\leq n,$}\\
L_i& \hbox{\small if $X\neq C(T,T'),C(1,T')$ and $(\dot R^i_{re})_{lg}\neq\emptyset; 1\leq i\leq n,$}\\
 E_i & \hbox{\small if $X\not=BC(T,T')$ and  $(\dot R^i_{re})_{ex}\neq\emptyset; 1\leq i\leq n,$ }\\
  S_2 & \hbox{\small if $X= C(1,T').$ }\\
\end{array}\right.$$
 Finally for types $B(1,T), B(T,1), $ (\ref{rel-s-alpha}) together with (\ref{equal}) implies  that $$F+2S_2=F+2S_1=S_{-\ep_0+\ep_t}-2\frac{(\ep_0,-\ep_0+\ep_t)}{(\ep_0,\ep_0)} S_{\ep_0}\sub S_{\ep_0+\ep_t}=F;\;\;\; (t\in T),$$ and  for types $BC(T,T')$  ($|T'|\geq2$) and $B(T,T'),$ (\ref{triple})($a$) together with   (\ref{equal}) and (\ref{*}) implies that  $$F+2S_1=F+2 S_2\sub F+L_2\sub F,$$ also for other types by (\ref{com}), we have $$F+2S_i\sub F+S_i+S_i\sub F+S_i\sub F;\; (1\leq i\leq n).$$
 This completes the proof of Claim 5.
Now we are ready  to complete the proof.

$(i)$ Assume that  {\small $X\neq A(\ell,\ell),BC(T,T'),C(T,T'),C(1,T').$} If {\small $X\neq B(T,T'),$ $B(1,T)$} or {\small$ B(T,1),$} by (\ref{888}), $S:=F=S_i$ ($1\leq i\leq n$) is a subgroup of $A^0$ and so $F+S\sub S.$ Now if {\small $X=  B(T,T'),B(1,T),B(T,1),$}  by (\ref{com}) and (\ref{equal}), $S:=S_1=S_2$ and $F+S\sub S.$ This together with (\ref{88}), (\ref{888}), (\ref{8}) and Proposition \ref{super-root} completes the proof.

$(ii)$ Let {\small $X=BC(1,T),BC(T,T')$ with $|T|,|T'|>1.$} Then $S:=S_1=S_2$ by (\ref{equal}) and so by (\ref{com}), $F+S\sub S.$ Also for $i\in\{1,2\},$ by (\ref{*}), $F+E_i\sub F$  and by (\ref{888}), $F=L_i$ if $(R^i_{re})_{lg}\neq\emptyset.$ Therefore we have $2F+E_i\sub E_i$ if   $(R^i_{re})_{lg}\neq\emptyset.$ Finally by
(\ref{triple}), we have $S_i+E_i\sub S_i$ and $E_i+4S_i\sub E_i.$
Now we are done using (\ref{8}), (\ref{88}), (\ref{***}) and proposition \ref{super-root}.

$(iii)$ Let {\small $X=C(1,T')$ with $|T'|>2.$}  Taking $S:=S_1,$ we have $F+S\sub F$ by (\ref{com}). Also by (\ref{888}), $F=S_2,$ so we are done using (\ref{triple}), (\ref{8}), (\ref{***}) and Proposition \ref{super-root}.

$(iv)$ Let {\small $X=C(T,T')$ with $|T|,|T'|>2.$} Using (\ref{888}), we have $F=S_1=S_2$ and so we are done using (\ref{triple}) together with (\ref{8}), (\ref{***}) and Proposition \ref{super-root}.
\qed

\end{document}